\documentclass[twoside,12pt]{article}
\usepackage{amsthm,amsmath,amsfonts,amscd,amssymb,dsfont}
\usepackage{latexsym,amsmath}
\usepackage[T1]{fontenc}
\usepackage[pagebackref=false,colorlinks,linkcolor=black,citecolor=magenta]{hyperref}
\usepackage{framed}
\usepackage{pgfplots}
\usepackage{tikz-cd}
\usepackage{mathrsfs}
\usepackage[all]{xy}
\usetikzlibrary{arrows}
\usepackage{hyperref}
\usepackage{mathrsfs}
\usepackage{upgreek}
\usepackage{blkarray}
\usetikzlibrary{arrows}
\usepackage{pgfplots}
\pgfplotsset{compat=1.7}
\usepackage{etoolbox,lastpage}
\newtheorem{theorem}{Theorem}[section]
\newtheorem{lemma}[theorem]{Lemma}
\newtheorem{proposition}[theorem]{Proposition}

\newtheorem{definition}[theorem]{Definition}
\newtheorem{remark}[theorem]{Remark}

\numberwithin{equation}{section}
\allowdisplaybreaks

\makeatother
\begin{document}
	
	\title{\large\bf Supertransversality and $\Pi$-symmetric supermanifolds}	
	\bigskip
	\author{{\small \bf Fatemeh Alikhani $^1$, Mehdi Ghorbani $^1$, Saad Varsaie $^1$\footnote{
				{ Corresponding author, E-mail address: \url{varsaie@iasbs.ac.ir}}\\$^1$ Institute for Advanced Studies in Basic Sciences, 444 prof. Yousef Sobouti Blvd, Zanjan 45137-66731, Iran\\}}~\\[5mm]}
	\date{}
	\maketitle	
	\begin{abstract}
		\noindent The main objective of this article is to extend the concept of transversality to supergeometry. Transversality has two important properties in the classical case, namely " stability" and " genericity", which we show in the following that in the category of smooth supermanifolds, supertransversality has stable property. By extending Sard's theorem to supergeometry, genericity property is proved. In the final section, we examine transversality in the category of
		$\Pi$-symmetric supermanifolds. The theory
		presented here is a step towards an extension of the concept of Euler-Poincaré characteristic to supermanifolds.
	\end{abstract}
	
	\bigskip
	\baselineskip=0.2cm
	\tableofcontents
	\bigskip
	\baselineskip=0.5cm
	\section{Introduction}\label{sec:introduction}
	Transversality is one of the fundamental concepts in differential topology. The concept of transversality was first introduced by the French mathematician René Frédéric Thom in 1950. Transversality appears in the context of cobordism theory \cite{15,16} and is a starting point for surgery theory \cite{17}. Transversality has other applications in topology \cite{18}, kinematics\cite{19}, physics \cite{20}, integrodifferential equation \cite{21}, economics \cite{22}, etc \cite{23,24}. Transversality is a concept that explains the intersection of a map and a submanifold, and the special case of transversality is the regular value theorem. In general, given a map
	$f: X \to Y$
	and a submanifold
	$Z \subset Y$,
	the goal is to investigate conditions under which the set
	$f^{-1}(Z)$
	is a submanifold of
	$X$.
	In this case, the map
	$f$
	is called transversal to the submanifold
	$Z$.
	Thus, transversality can be considered as a property of maps. In differential topology, Thom's transversality theorem is a major result that says transversality is a stable and generic property \cite{3}.
	In this article, we first define the notion of supertransversality for morphisms between supermanifolds and then prove part of Thom's theorem in the category of supermanifolds. One of the important features of transversality is "stability". By a stable property of maps we mean it is preserved under small deformation. In this article, we attempt to prove some stable properties for morphisms between supermanifolds while defining the notion of superhomotopy. Among them, we have shown that supertransversality is a stable property, which we discuss in section four. Another important property of transversality is "genericity", meaning that any arbitrary map can be transformed into a map that is transversal to a given submanifold with a small change. This important property will be 
	proved in
	section 6, by using a super generalization of Sard's theorem. Indeed, we have extended Sard's theorem in supergeometry. In conclusion, we studied transversality and its associated properties within the category of $\Pi$-symmetric supermanifolds
	\cite{1,2,4,5}.\\
	The precise definition and interpretation of topological invariants such as the Euler-Poincaré characteristic has been a longstanding challenge in the general theory of supermanifolds. Establishing a rigorous framework for intersection theory and transversality in this supersymmetric geometric setting is fundamental for both mathematics and theoretical physics.
	As mentioned above, in this paper, the concept of supertransversality was introduced and rigorously developed in the broader category of supermanifolds so this theory
	provides essential tools to analyze generic intersections and topological invariants. The theory has been further refined and specialized to the subclass of orientable $\Pi$-symmetric supermanifolds\cite{25}.
	In this context, an intrinsic definition of the Euler-Poincaré characteristic pair has been formulated specifically for orientable $\Pi$-symmetric supermanifolds. This new topological invariant not only addresses a previously open problem in supergeometry but also paves the way for applications in theoretical physics, particularly in the study of supermoduli spaces, the Batalin-Vilkovisky formalism, topological quantum field theory, and supergravity. Inspired by pioneering works such as Witten’s 1988 development of topological quantum field theory, this advancement provides foundational insights for further development of topological structures in quantum field and string theories\cite{26}.
	\section{Preliminaries}
	In this section, some necessary concepts and results are introduced as follows:\\
	\textbf{Ringed space}.
	Let $X$ be a topological space and $R$ be a sheaf of rings i.e. for every open subset $U$ of $X$, $R(U)$ is a commutative ring. Then $(X, R)$ is called a ringed space.  The elements of $R(U)$ are called the \textit{sections} of $R$ over $U$. We define the \textit{stalk} $R_x$ of $R$ at $x\in X$ to be the direct limit of the rings $R(U)$ for all open sets $U$ containing $x$ via the restriction maps. A ringed space is said to be a \textit{space} if the stalks are local, meaning that each stalk $R_x$ has a unique maximal ideal denoted by $\mathit{m}_x$.
	\\
	\textbf{Products}.
	The category of supermanifolds admits products. Let $X_{i}\, (1 \leq i\leq n)$ be spaces in the category. A ringed space $X$ together with projection maps $P_{i}$ : $X\rightarrow X_{i}$ is called a product of the $X_{i},$
	$$
	X=X_{1}\times\cdots\times X_{n},
	$$
	if the following is satisfied: for any ringed space $\mathrm{Y}$, the map
	\begin{equation}\label{2.1..}
		f\mapsto(P_{1}\circ f,\ \ldots,\ P_{n}\circ f)
	\end{equation}
	from $\mathrm{H}\mathrm{o}\mathrm{m}(\mathrm{Y},\ X)$ to $\prod_{i}\mathrm{H}\mathrm{o}\mathrm{m}(Y,\ X_{i})$ is a bijection. In other words, every morphism $f$ from $\mathrm{Y}$ to $X$ are identified with $n$-tuples $(f_{1},\ \ldots,\, f_{n})$ of morphisms $f_{i}(\mathrm{Y}\rightarrow X_{i})$. This product is unique up to isomorphism \cite{1}.\\
	\textbf{Superdomain}.
	Let $U$ be an open set in $\mathbb{R}^p$, and let $\mathcal{C}^{\infty \, p|q}$ be the sheaf of smooth supercommutative $\mathbb Z_2$-graded rings on $U$ that assigns to each open $V \subseteq U$ the ring	$ C^\infty (V) \,\, [ \theta^1, ..., \theta^q ]$ i.e. the ring of polynomials in the variables $ \theta^1 ,..., \theta^q$, satisfying the relations $ \theta^i \theta^j= - \theta^j \theta^i$ and $\theta^{i^2}=0$ and coefficients in $C^{\infty}(V)$ the ring of smooth real valued functions on $V$. Then the pair $(U, \mathcal{C}^{\infty \, p|q} )$ is called a superdomain, denoted by $U^{p|q}$.
	\begin{remark} Each element of $\mathcal{C}^{\infty p|q}(V)$ can be written as 
		\begin{align*}
			\sum\limits_{I \subset\{1,\dots,q\}} f_I\theta^I,
		\end{align*}
		where the $f_I \in C^{\infty}(V)$ and $\theta^{I}=\theta^{i_1} \cdots \theta^{i_r}$ if $I=\{i_1,\dots , i_r\}$ such that $\ i_1< \dots < i_r$.
	\end{remark}
	\textbf{Supermanifold}.
	A supermanifold of dimension $m|n$ is a ringed space 
	$X=(\tilde{X},\mathcal{O}_X)$ that is locally isomorphic to $U^{m|n}$ and $\tilde{X}$ is a second countable
	and Hausdorff topological space. The superdomains $\mathbb{R}^{m|n}$ and $U^{m|n}$ are special examples of supermanifolds.\\
	\textbf{Morphism}.
	A morphism $\varphi$ from the supermanifold $X=(\tilde X,\mathcal O_X)$ to the supermanifold $Y=(\tilde Y,\mathcal O_Y)$ consists of a pair $(\tilde \varphi, \varphi^*)$, where $\tilde \varphi:\tilde X\to \tilde Y$ is a continues function and
	$\varphi^*: \mathcal O_Y \to \varphi_*(\mathcal O_X)$
	is a homomorphism between sheaves of supercommutative $\mathbb Z_2$- graded rings. In other words, for each open  $\tilde V \subset \tilde Y$, there is a homomorphism
	$$\varphi^*_{\tilde V}: \mathcal O_Y (\tilde V) \to \mathcal O_{X} ( \tilde \varphi^{-1} (\tilde V))$$
	such that for every open subsets $\tilde U$ and $\tilde V$ with $\tilde U\subset \tilde V$ the following diagram commutes
	\begin{equation}
		\begin{CD}
			\mathcal{O}_Y(\tilde V)@>\varphi^*_{\tilde V}>>\mathcal{O}_X(\tilde \varphi^{-1} (\tilde V))@.\qquad \\
			@ Vr^Y_{\tilde U\tilde V}VV @ VV r^X_{\tilde U\tilde V} V @.\\
			\mathcal{O}_Y(\tilde U)@>>\varphi^*_{\tilde U}>\mathcal{O}_X(\tilde \varphi^{-1} (\tilde U))@.
		\end{CD}.
	\end{equation}
	This collection of homomorphisms can be thought of as a set of constraints that relate the structure sheaf of $X$ to the structure sheaf of $Y$. It is worth noting that by these constraints it can be shown that $\tilde \varphi$ is a smooth map if $\tilde{X}$ and $\tilde{Y}$ are equipped with their differential structures. Indeed, One can show that there exists a unique differential structure on $\tilde X$ under which $\tilde X$ is a smooth manifold embedded in $X$.\\
	We have presented other concepts of supergeometry in the appendix to concentrate on the main topics of the article. For more details on the geometry of supermanifolds, the reader is suggested to
	refer to books \cite{1, 6, 7, 8, 9, 10}.
	\section{Supertransversality}\label{3}
	In this section, we first recall the transversality in classical differential topology and then introduce supertransversality.
	\begin{definition}(\textbf{Transversality})
		Let $f:X \to Y$ and $Z$ be a submanifold of $Y$. The map $f$ is said to be transversal to the submanifold $Z$, denoted by $f \pitchfork Z$, if the equation 
		$$	Im\,df_x+T_yZ=T_yY$$ 
		holds for all points $x$ in the preimage $f^{-1}(Z)$ and $y=f(x).$
	\end{definition}
	\begin{theorem} If the smooth map $f:X \to Y$ is transversal to a submanifold $Z \subset Y$, then the preimage $f^{-1}(Z)$ is a submanifold of $X$.\\ Moreover, the codimension of $f^{-1}(Z)$ in $X$ equals the codimension of $Z$ in $Y$\cite{3}.
	\end{theorem}
	In the next definition, we introduce the concept of supertransversality.
	\begin{definition}(\textbf{Supertransversality})
		Let 
		$ \psi=(\tilde \psi,\psi^*):X=(\tilde X,\mathcal{O}_X)\longrightarrow Y=(\tilde Y,\mathcal{O}_Y) $
		be a morphism between supermanifolds, and
		$ Z=(\tilde Z,\mathcal{O}_Z) $
		be a subsupermanifold of $ Y $. The morphism $\psi$ is said to be supertransversal to $Z$, if for each point $x$ in $\tilde \psi^{-1}(\tilde Z)$ and $ y=\tilde{\psi}(x)$ the equation
		$$ 	Im\,d\psi_x+T_yZ=T_yY $$ holds true. Then as in the classical way, denoted by $\psi \pitchfork Z$.
	\end{definition}
	\begin{remark}
		Let $W$ be a subsupermanifold of $X$. The codimension of $W$ in $X$, denoted by $codim_X^ W$,  is said to be $r|s$ if $dim{X}=m+r|n+s$ and $dim{W}=m|n$.
	\end{remark}
	\begin{theorem}(\textbf{Transversality in supergeometry})\label{3.5.}
		Let
		$ \psi:X^{m+r|n+s}\longrightarrow Y^{p+r|q+s} $
		be a morphism between supermanifolds, and
		$ Z^{p|q} $ is a subsupermanifold of $ Y $ such that 
		$ \psi\pitchfork Z $. Then there exists a subsupermanifold $ W^{m|n} $
		of $ X $ (henceforth it will be denoted  by $\psi^{-1}(Z)$ for convenience)
		and a morphism
		$\hat{\psi}$
		such that the following diagram is commutative:
		\begin{equation} \label{1}
			\begin{CD}
				X^{m+r|n+s}@>\psi>>Y‌^{p+r|q+s}\\
				@ A jA A @ AAiA\\
				W^{m|n}@ >>\hat{\psi}> Z^{p|q}
			\end{CD},
		\end{equation}
		where $i$ and $j$ are immersions with $\tilde{i}$ and $\tilde{j}$ are inclusions.
		Moreover,
		$$codim^{W}_{X}=codim^{Z}_{Y}=r|s.$$
	\end{theorem}
	\begin{proof}
		As a first step, for each $y\in\tilde{Z}$, we construct a submersion
		$$ \Psi_{y}:V_{y}^{m+r|n+s}\longrightarrow \mathbb R^{r|s},$$
		where $V_y^{m+r|n+s}$ is an open subsupermanifold of $X$, with
		$\tilde\Psi_y^{-1}(0)= \tilde{\psi}^{-1}(\tilde{Z})\cap\tilde{V}_{y}$.
		Since $Z^{p|q}$ is subsupermanifold of $Y$, according to \ref{2.10.}, for each point $y\in \tilde{Z}$, there are isomorphisms of supermanifolds $\varphi$ and $\varphi'$ from an open neighborhood of $y$, say $U$, in $Y$ and $U\cap Z$ in $Z$ two open  neighborhoods of origin in $\mathbb R^{p+r|q+s}$ and $\mathbb R^{p|q}$ respectively, such that the following diagram commutes: 
		\begin{equation} \label{(2)}
			\begin{CD} 
				{U\subset Y^{p+r|q+s}}@>\varphi>>\bar{U}\subset\mathbb{R}‌^{p+r|q+s}\\
				@ A iA A @ AA\bar{i}A\\
				{U\cap Z}\subset Z^{p|q}@ >>{\varphi^{'}}>\overline{U\cap Z}\subset \mathbb{R}^{p|q}
			\end{CD}.
		\end{equation}
		By $\bar{U}$ and $\overline{U\cap Z}$ we mean standard superdomains corresponding to $U$ and ${U\cap Z}$ respectively.
		Now we consider standard coordinate system  $(t_1,\dots,t_{p+r};e_1,\dots,e_{q+s})$ for $\mathbb {R}^{p+r|q+s}$. Then one has \ref{2.10.}
		$$\bar{i}^*: \mathcal{O}_{\mathbb{R}^{p+r|q+s}}\to  \mathcal{O}_{\mathbb{R}^{p|q}} $$
		\begin{equation*}
			\begin{cases}
				t_i \mapsto t_i \qquad 1\leq i \leq p &\\
				t_i \mapsto 0  \qquad p+1 \leq i \leq p+r &\\
				e_j \mapsto e_j  \qquad 1 \leq j \leq q &\\
				e_j \mapsto 0 \qquad q+1 \leq j \leq q+s
			\end{cases}.
		\end{equation*}
		Set $G_y= Pr\circ \varphi$, where $Pr$ is a projection as follows:
		$$ Pr: \bar{U}\subset{\mathbb{R}^{p+r|q+s}}\to  \bar{U'}\subset{\mathbb{R}^{r|s}} $$
		\begin{equation*}
			\begin{cases}
				t_i \mapsto 0 \qquad 1\leq i \leq p &\\
				t_i \mapsto t_{i_p}  \qquad p+1 \leq i \leq p+r &\\
				e_j \mapsto 0  \qquad 1 \leq j \leq q &\\
				e_j \mapsto e_{j_q} \qquad q+1 \leq j \leq q+s
			\end{cases}.
		\end{equation*}
		It can be easily seen that $G_y$ is a submersion.
		Now let
		$ \Psi_{y}:V_{y}\subset X^{m+r|n+s}\longrightarrow \mathbb R^{r|s}$
		be the morphism $\Psi_{y}=G_y\circ \psi$. Since by assumption $\psi \pitchfork Z$ and $G_y$ is submersion, it follows that $\Psi_{y}$ is a submersion.
		Now we consider $\Psi_{y}^{-1}(0)$,
		$$\tilde{\Psi_{y}}^{-1}(0)=(\tilde{G}_y\circ\tilde{\psi})^{-1}(0)=\tilde{\psi}^{-1}(\tilde{U}_{y}\cap\tilde{Z})\subset\tilde{W}\subset\tilde{X}.$$
		For convenience, we set $\tilde{W}_y:=\tilde{\psi}^{-1}(\tilde{U}_y\cap\tilde{Z}).$
		The theorem of level sets of submersions in supergeometry(c.f. section 4.7 \cite{1}) shows that $\Psi_y^{-1}(0)$ is subsupermanifold of $X$. And its underlying topological space is $\tilde{W}_y$. Set $W_y:=\Psi_y^{-1}(0).$ 
		One has $\tilde W=\bigcup \tilde \psi^{-1}(\tilde U_y \cup \tilde Z) =\bigcup_y \tilde W_y$. We show that $W$ can be obtained by gluing $W_y$'s. By submersion level set theorem in supergeometry, for each $y$, $W_y$ is a subsupermanifold of $X$. It is sufficient to introduce morphisms that gluing $W_y$'s along them gives $W$.
		\begin{equation} \label{(3)}
			\begin{CD}
				\mathcal{O}_X(\tilde V)@>id>>\mathcal{O}_X(\tilde V)@.\qquad \\
				@ Vi^{*}_yVV @ VV {i^{*}_{y^{'}}} V @.\\
				\mathcal{O}_{W_{y}}(\tilde W_{y}\cap \tilde V)@>>f^{*}_{y{y^{'}},x}>\mathcal{O}_{W_{y'}}(\tilde W_{y'}\cap \tilde V)@.
			\end{CD}.
		\end{equation}
		Let $x \in \tilde W_y \cap \tilde W_{y'}$. Since $W_y$ and $W_{y'}$ are subsupermanifolds of $X$, the theorem \ref{1.6.} shows that for an open subsupermanifold of $X$ say $V$,  there are sets of generators $(f_i, g_j)$, $(h_i, k_j)$ for $\mathcal O_X(\tilde{V})$ and $(f'_t, g'_l)$, $(h'_t, k'_l)$ for $\mathcal O_{W_y}(\tilde{W_y}\cap{\tilde V})$ and $\mathcal O_{W_{y'}}(\tilde{W_{y'}}\cap{\tilde V})$ respectively, such that with respect to them the natural embeddings $i_y: W_y\to X$ and $i_{y'}: W_{y'}\to X$ are as follows:
		\begin{center}
			\begin{tabular}{cc}
				$i_{y}^*(f_i)=	\begin{cases}
					f_i^{'} \qquad 1 \leq i \leq m \\\
					0  \qquad m+1 \leq i \leq m+r
				\end{cases}$& 
				\text{and} $\qquad	i_{y}^*(g_j)=	\begin{cases}
					g_j^{'} \qquad 1 \leq j \leq n\\\
					0  \qquad n+1 \leq j \leq n+s
				\end{cases}$
			\end{tabular}
		\end{center}
		Now one may define a morphism say $f_{yy',x}=(id_{\tilde{W}_y\cap \tilde{V}}, f^*_{yy',x})$
		such that commutes the diagram \ref{(3)}. The family of  morphisms $\{f_{yy',x}, x\in \tilde{W}_y\cap\tilde{W}_{y'}\}$, defines a  morphism  $f_{yy'}:W_y\to W_{y'}$. Gluing $W_y$ through morphisms $f_{yy'}$ give the subsupermanifold $W$.
		For the existence of $\hat{\psi}$, it is sufficient to note that it exists on $W_y$ for each $y$.
	\end{proof}
	In the next section, we first generalize the concept of homotopy between morphisms in supergeometry and then use it to prove stability for some properties, particularly supertransversality.
	\section{Stability}
	Many properties of maps remain invariant under smooth deformation. A smooth map is said to be smoothly homotopic to another if they can be smoothly deformed into one another using a family of smooth maps. The precise formulation of this concept is one of the most important topics in topology and is called \textit{homotopy}. Homotopy is an equivalence relation on smooth maps, and the equivalence class to which a map belongs it is called its homotopy class.
	
	Certain properties of maps are physically meaningful if remain invariant under small deformations of the map. Such properties are called \textit{stable properties}. A collection of maps possessing a special stable property can be referred to as a class of stable maps.
	In particular, a stable property is one that, given $f_0:X \to Y$ possessing that property, and a homotopy $f_t:X \to Y$ with $f_0$ as its initial map, then for some $\epsilon > 0$, all $f_t$ with $t < \epsilon$ also possess that property. In classical differential topology, it is established that the collection of local diffeomorphisms, immersions, submersions, maps transversal to a closed submanifold, embeddings, and diffeomorphisms of a compact manifold form a class of stable maps \cite{3}.
	In this section, we aim to define \textit{superhomotopy} and prove some stable properties for morphisms between supermanifolds.
	\begin{definition}\label{4.1.}
		Suppose
		$\psi_0: X^{m|n} \to Y^{p|q}$
		and
		$\psi_1:X^{m|n} \to Y^{p|q}$
		are morphisms between supermanifolds. We say
		$\psi_0$
		and
		$\psi_1$
		are homotopic, and denote it by
		$\psi_0 \stackrel{\Psi}{\sim} \psi_1$ if there exists a morphism
		$$\Psi:X \times \mathbb{R}^{1|0} \to Y$$
		such that
		$\Psi \circ i_0=\psi_0$
		and
		$\Psi \circ i_1=\psi_1$,
		where
		$ i_0\, \text{and}\,\, i_1: X^{m|n} \to X^{m|n} \times \mathbb{R}^{1|0}$
		are morphisms defined as follows:
		$$	\begin{tabular}{cc}
			$\begin{cases}
				\tilde i_0: \tilde X \to \tilde X \times \mathbb{R}^{1}\\
				\,\,\,\,\,\,\,\,\,\,\,\, x \mapsto (x,0)\\
			\end{cases}$
			&$\begin{cases}
				\tilde i_1: \tilde X \to \tilde X \times \mathbb{R}^{1}\\
				\,\,\,\,\,\,\,\,\,\,\,\, x \mapsto (x,1)\\
			\end{cases}$
		\end{tabular}$$
		also, the morphisms $i_0^*$ and $i_1^*$ are defined locally in terms of the coordinate system \\
		$(x_1, \ldots, x_m, t; e_1, \ldots, e_n)$, where $t$ is a coordinate to the supermanifold $\mathbb{R}^{1|0}$, as follows:
		\begin{center}
			\begin{tabular}{cc}
				$\begin{cases}
					i_0^*:f(x_1,\dots,x_m,t ) \mapsto f(x_1,\dots,x_m,0)\\
					i_0^*:e_j \mapsto e_j
				\end{cases}$&	$\begin{cases}
					i_1^*:f(x_1,\dots,x_m,t ) \mapsto f(x_1,\dots,x_m,1)\\
					i_1^*:e_j \mapsto e_j
				\end{cases}.$
			\end{tabular}
		\end{center}
	\end{definition}
	\begin{remark}
		The category of supermanifolds has the property of categorical products (c.f. \ref{2.1..}). Therefore, 
		$X \times \mathbb{R}^{1|0} = (\tilde{X} \times \mathbb{R}, \mathcal{O}_{X \times \mathbb{R}^{1|0}})$ 
		means that it is the supermanifold obtained from the product of $X$ and $\mathbb{R}^{1|0}$ with each other. Moreover, if the coordinate system $(x_1, \dots, x_m; e_1, \dots, e_n)$ is defined on $X$ and ${t}$ is the coordinate system defined on $\mathbb{R}^{1|0}$, then $(x_1, \dots, x_m, t; e_1, \dots, e_n)$ is a coordinate system on $X \times \mathbb{R}^{1|0}$.
	\end{remark}
	\begin{definition}
		We say that a morphism $\psi_t$ is a superhomotopy of $\psi_0:X^{m|n} \to Y^{p|q}$ if there exists a morphism $\Psi: X \times \mathbb{R}^{1|0} \to Y$ such that $\Psi \circ i_0 = \psi_0$ and $\Psi \circ i_t = \psi_t$, where $i_t: X^{m|n} \to X^{m|n} \times \mathbb{R}^{1|0}$ is a morphism similar to the morphisms $i_0$ and $i_1$ in definition \ref{4.1.}.
	\end{definition}
	In the following theorem, we prove that the property of submersion in the category of supermanifolds is stable under certain conditions.
	\begin{theorem}
		Let
		$ \psi_0 : X^{m|n} \to Y^{p|q}$
		be a submersion. If
		$\tilde{X}$
		is compact and
		$ \psi_t $
		is a superhomotopy of
		$\psi_0$
		, then there exists an
		$\varepsilon >0$
		such that for every
		$t$
		smaller than
		$\varepsilon$,
		$ \psi_t:X^{m|n} \to Y^{p|q}$
		is also a submersion.
	\end{theorem}
	\begin{proof}
		Consider an arbitrary point $(x_1,t_1)$ in $\tilde X \times \mathbb R^1$. Assume that $(d\psi_{t_1})_{x_1}$ is surjective. In this case, consider open neighborhoods $\tilde U_{x_1}\subset \tilde X$ with coordinates $(x^1,...,x^m; e^1,...,e^n)$ around $x_1$, and $\tilde V_{y_1}\subset \tilde Y$ with coordinates $(y_1,...,y_p;\theta_1,...,\theta_q)$ around $y_1=\tilde \psi(x_1)$ such that $\tilde{\psi}(\tilde U_{x_1})\subset\tilde V_{y_1}$. If
		$$\psi_{t_1}^*:\mathcal O_{Y,y_1} \to \mathcal O_{X,x_1}$$
		$$\begin{cases}
			y_i \mapsto f_i\\
			\theta_j \mapsto g_j
		\end{cases}$$
		then for $	1\leq a \leq p, 	1\leq b \leq q, 1\leq i \leq m, 1\leq j \leq n$ we have
		\begin{align}
			[ d \psi_{t_1}]_{x_1}  &= \left(\begin{array}{cc}
				\dfrac{\partial f_a}{\partial x^i}(t_1)&-\dfrac{\partial f_a}{\partial e^j}(t_1)\\
				\\
				\dfrac{\partial g_b}{\partial x^{i}}(t_1)&   \,\, \dfrac{\partial g_b}{\partial e^{j}}(t_1)\\
			\end{array} \right)
		\end{align}	
		Let $\mathcal{J}$ be ideal of nilpotent elements of $\mathcal{O}$ and let $[d\psi_{t_1}]^{\sim}$ be a matrix over $\mathcal{O}/\mathcal{J}$ obtained  by applying the map $\mathcal{O}\to \mathcal{O}/\mathcal{J}$ to each entry of $[d\psi_{t_1}]$. Thus
		\begin{align}
			[ d \psi_{t_1}]_{x_1}^{\sim}  &= \left(\begin{array}{cc}
				(  \dfrac{\partial f_a}{\partial x^i})^{\sim}(t_1)&0\\
				0& (  \dfrac{\partial g_b}{\partial e^{j}})^{\sim}(t_1)\\
			\end{array} \right).
		\end{align}
		Due to the submersion property of $\psi_{t_1}$, $ ( \dfrac{\partial f_a}{\partial x^i})^{\sim}(t_1)$ has a non-singular $p\times p$ submatrix say $(I_p(t_1))$ and $ (\dfrac{\partial g_{b}}{\partial e^{j}})^{\sim}(t_1)$ also has a non-singular $q\times q$ submatrix $(J_q(t_1))$.
		Since the
		$\psi_{t_1}=\Psi \circ i_{t_1}$
		and
		$i_{t_1}^* \circ \Psi^*(y_i)=f_i$
		, it follows that
		$f_i$
		and its partial derivatives are smooth functions of
		$t_1$. So, the determinants
		$I_p(t_1)$
		and
		$J_q(t_1)$
		are smooth functions of $t_1$. As a result, there exists a neighborhood $\tilde U_{x_1}\times I_{t_1}$ around the point $(x_1,t_1)$, such that for $t \in I_{t_1}=(t_1-\varepsilon,t_1+\varepsilon)$, $\psi_t$ is a submersion.
		Now according to the assumption
		$\psi_{t_1=0}$
		is a submersion. So, for
		$x\in \tilde X$
		there exists a neighborhood
		$\tilde U_x \times (-\varepsilon_x, \varepsilon_x)$
		around
		$(x,0)$
		such that
		$\psi_t$
		is a submersion for
		$t\in(-\varepsilon_x, \varepsilon_x)$
		. Therefore,
		$\mathcal U=\bigcup(\tilde U_x \times (-\varepsilon_x, \varepsilon_x))$
		is a neighborhood for
		$\tilde X \times \{0\}$
		and since
		$\tilde X$ is compact, there exists
		$\varepsilon$
		such that
		$\tilde X \times [0, \varepsilon)\subset \mathcal U$. Thus, for all
		$t < \varepsilon$,
		$\psi_t$
		is a submersion.
	\end{proof}
	In the following theorem, we demonstrate that transversality is a stable property in supergeometry.
	\begin{theorem}
		Assume a morphism $\psi_0 : X^{m+r|n+s} \to Y^{p+r|q+s}$ is supertransversal to subsupermanifold $Z^{p|q}$ of $Y$. If
		$\tilde X$
		is compact
		and
		$\tilde Z$
		is a closed subsupermanifold of
		$\tilde Y$
		and
		$\psi_t$
		be a superhomotopy of
		$\psi_0$, then there exists
		$\varepsilon >0$
		such that for every
		$t$
		smaller than
		$\varepsilon$,
		$ \psi_t$
		also be a supertransversal to
		$Z$
		.
	\end{theorem}
	\begin{proof}
		The idea of the proof is as follows:
		Since
		$\psi_0 \pitchfork Z$
		, so for every
		$ y\in \tilde{Z}$
		, there exists a morphism
		$$ G_y: U_y^{p+r|q+s}\to \mathbb R^{r|s}$$
		such that 
		$ G_y\circ \psi_0$
		is a submersion. Here
		$\tilde U_y$
		is an open neighborhood of
		$y$
		in
		$ \tilde{Y}$.
		So, if
		$G_y\circ \psi_0$
		be a submersion then there exists $\varepsilon_y>0$ such that for each $t < \varepsilon_y$, $G_y\circ \psi_t$ is also a submersion.
		Obviously, $\{\tilde U_y\}$
		is an open covering for
		$\tilde Z$. On the other hand,
		$\tilde X$ is compact so
		$\tilde{\psi_0}(\tilde X)$
		is compact. As a result, the intersection of the latter set with $\tilde Z$ is also compact, so it has a finite subcover such as $\left\lbrace \tilde U_{y_i}\right\rbrace _{i=1,...,l}$.  Put
		$ \varepsilon = min\{\varepsilon_1, ...,\varepsilon_l\}$. Thus, for each $t< \varepsilon$, $G_y\circ \psi_t$ is a submersion, and equivalently, $\psi_t$ is transversal with respect to $Z$.
	\end{proof}
	In the subsequent sections, we show that a generalization of Sard's theorem is held in supergeometry and then use it to establish the genericity property of transversality in supermanifolds.

	\section{Sard's theorem in Supergeometry}
	In this section, we first introduce a vector bundle associated with the structure
	sheaf of a supermanifold. We then define the notions of regular and critical values and present a generalization of Sard's theorem within the context of supergeometry.

	Let
	$M=(\tilde M,\mathcal O_M)$
	be a supermanifold and let $\mathcal J$ be the subsheaf of ideals generated
	by odd elements. Therefore
	$\left( \dfrac{\mathcal J}{\mathcal J^2}\right) $ is a locally free
	$\mathcal{O/J}$- module.
	that
	$\mathcal J$
	is an ideal maximal ring of
	$\mathcal O_M$.
	\begin{lemma}\label{3.1}
		Let
		$M^{m|n}=(\tilde M,\mathcal O_M)$
		be a supermanifold. Then there exists a vector bundle
		associated with $ \dfrac{\mathcal J}{\mathcal J^2}$.
	\end{lemma}
	\begin{proof}
		At any point
		$x$
		we consider
		the stalk
		$\mathcal{J}_x/\mathcal{J}^2_x$.
		Suppose
		$\xi_i$ for $1\leq i\leq n$, be odd coordinates in $U_{\alpha}$ a neighborhood of
		$x$.
		We consider vector space $E_x=\left\langle [\xi_i]_x+\mathcal J^2_x\right\rangle_\mathbb R $
		and put\\
		$E=\bigsqcup_{x \in \tilde M} E_x$. It can be shown
		$$\Xi=(E, \pi, \tilde M, \mathbb R^n)$$
		is a vector bundle such that its sheaf of sections is isomorphic with
		$\dfrac{\mathcal J}{\mathcal J^2}$.\\
		In the following, we will consider the local charts of a vector bundle
		$\Xi$.
		Set:
		\begin{align*}
			&h_\alpha: U_\alpha \times \mathbb R^n \to \pi^{-1}(U_\alpha)\\
			& (x,(v_1,...,v_n)) \mapsto \sum_{i=1}^{n}v_i([\xi_i]_x+\mathcal J^2_x)
		\end{align*}
		and
		\begin{align*}
			&h_\alpha^{-1}: \pi^{-1}(U_\alpha) \to U_\alpha \times \mathbb R^n \\
			& \sum_{i=1}^{n}\alpha_i([\xi_i]_x+\mathcal J^2_x) \mapsto (x,(\alpha_1,...,\alpha_n)).
		\end{align*}
		At any point
		$x \in U_\alpha \cap U_\beta$ we consider local coordinates
		$(x_1,...,x_m;\xi_1,...,\xi_n)$
		for
		$U_\alpha$
		and local coordinates $(y_1,...,y_m;\eta_1,...,\eta_n)$
		for
		$U_\beta$.
		So one has
		$\xi_i= \sum_{J\subset{1,...,n}}f_{iJ}\eta_J$.
		Therefore
		$$[\xi_i]_x +\mathcal{J}_x^2=\sum_{j=1,..., n}
		[f_{ij}]_x[\eta_j]_x + \mathcal{J}^2_x.$$
		The transition map is as follows:
		\begin{align*}
			h_{\beta}^{-1}\circ h_{\alpha} : & U_{\alpha\beta}\times \mathbb R^n \to U_{\alpha\beta}\times \mathbb R^n\\
			&(x,V)=(x,(v_1,...,v_n)) \mapsto \sum_{i} v_i([\xi_i]_x+\mathcal J^2_x)\\
			&\,\,\,\,\,\,\qquad\qquad\,\,\,\, \mapsto \sum_{i} v_i\left( (\sum_{j=1}^{n}	[f_{ij}]_x[\eta_j]_x)+\mathcal J^2_x\right)\\
			&\,\,\,\,\,\,\qquad\qquad\,\,\,\, =\sum_{j}\left(\sum_{i}  v_i 	[f_{ij}]_x\right) \left( [\eta_j]_x+\mathcal J^2_x\right)\\
			&\,\,\,\,\,\,\qquad\qquad\,\,\,\, \mapsto\left(x,\left(\sum_{i}  v_i	[f_{i1}]_x,...,\sum_{i}  v_i	[f_{in}]_x \right)  \right)\\
			&\,\,\,\,\,\,\qquad\qquad\,\,\,\, =\left( x,V(	[f_{ij}]_x)\right) 
		\end{align*}
		therefore
		$\Xi$ is a vector bundle.
	\end{proof}
	In what follows, we show that Sard's theorem holds in supergeometry. To this end, let's start with introducing some necessary notations and propositions., the second component of the intersection pair can be calculated as a classical term. \\
	In the next proposition, we show that the dual sheaf of the quotient sheaf
	$\left( \dfrac{\mathcal J}{\mathcal J^2}\right) _X$
	at the point
	$x$
	is isomorphic  with the odd subspace of the tangent space
	$T_xX$
	. For this purpose, we introduce some useful morphisms. \\
	Consider the morphism
	$\psi: X \to Y$.Then
	$$\psi^*: \mathcal{O}_{Y} \to \mathcal O_{X}.$$
	this morphism induces the following bundle morphism:
	$${\bar{{\psi^*}}}:\Xi_{Y}\to \Xi_{X}$$
	where $\Xi_{Y}$ is vector bundle associated to $(\mathcal{J}/\mathcal{J}^2)_{Y}$. On the dual bundles, we have the following bundle morphism:
	$$(\bar{\psi^*})^*: \left( {\dfrac{\mathcal J}{\mathcal J^2}}\right)^* _{X }\to \left( {\dfrac{\mathcal J}{\mathcal J^2}}\right)^*_{Y}$$ 
	Now consider the morphism 
	\begin{align}
		&i_x:\left( \left({\dfrac{\mathcal J}{\mathcal J^2}}\right) _X \right)_x ^* \to ({T_x}X)_1\\
		&\qquad \qquad \qquad \quad \,\varphi \mapsto \tilde {\varphi}\label{15}
	\end{align}
	where the derivation $\tilde{\phi}$ is defined as follows:
	\begin{equation}\label{16}
		\tilde{\varphi}(f):=\varphi((f-\tilde{f})+\mathcal J^2_X).
	\end{equation} Morphisms
	$i_z$
	and
	$i_y$
	are defined similarly.
	\begin{proposition}
		Morphisms
		$i_x$
		,
		$i_z$
		, and
		$i_y$
		defined above are isomorphisms.
	\end{proposition}
	\begin{proof}
		Every derivation
		$v \in ({T_x}X)_1$
		clearly defines a linear function
		$\varphi_v$
		on
		$\mathcal J_X$
		that is zero on
		$\mathcal J^2_X$.
		Conversely, for every
		$\varphi \in \left(\left(\dfrac{\mathcal J}{\mathcal J^2} \right)_{X}\right) _{x}^*$
		we have a tangent vector
		$v_\varphi$
		at
		$x$
		which is defined as
		$$v_\varphi(f)=\varphi\left( (f-\tilde f)+\mathcal J^2_X\right).$$
		The linearity condition for
		$v_\varphi$
		is clear. To check the Leibniz's property we have
		\begin{align*}
			v_\varphi(fg)=&\varphi\left( (fg-\tilde f\tilde g)+\mathcal J^2_X\right)\\
			=&\varphi\left( \left( (f-\tilde f)(g-\tilde g)+\tilde f(g -\tilde g)+(f-\tilde f)\tilde g\right) +\mathcal J^2_X\right)\\
			=&\tilde fv_\varphi(g)+(-1)^{p(f)p(g)}\tilde g v_\varphi(f).
		\end{align*}
		Therefore
		$v_\varphi$
		is a derivation. Also
		$v_{\varphi_v}$
		is
		equal to
		$v$.
	\end{proof}
	\begin{proposition} \label{5.9.0}
		Suppose
		$\psi:X \to Y$
		is a morphism between supermanifolds. Then for the bundle morphism
		\begin{align*}
			&(\bar{\psi^*})^*:\left(\left( \dfrac{\mathcal J}{\mathcal J^2} \right)_X\right)_x ^* \to \left(\left( \dfrac{\mathcal J}{\mathcal J^2} \right)_Y \right)_y ^*
		\end{align*}
		the following diagram is commutative.
		\begin{equation} \label{5.8.0}
			\begin{CD}
				({T_x}X)_1 @>d{\psi_x}>> ({T_y}Y)_1 \\
				@A{i_x}A{\rotatebox{90}{$\cong$}}A @A{\rotatebox{90}{$\cong$}}A{i_y}A \\
				\left( \left( \dfrac{\mathcal J}{\mathcal J^2} \right) _X \right)_x^* @>>(\bar{\psi^*})^*> \left( \left( \dfrac{\mathcal J}{\mathcal J^2} \right) _Y \right)_y^*
			\end{CD}
		\end{equation}
	\end{proposition}
	\begin{proof}
		For this purpose, we show
		$(d{\psi})_1 \circ i_x= i_y \circ (\bar{\psi^*})^*$.
		For
		$g \in \mathcal{O}_Y$
		on the left-hand side, we have:
		$$(d{\psi})_1 \circ i_x(\varphi)(g)=(d{\psi})_1 \circ {v_\varphi}(g)={v_\varphi} \circ \psi^*(g)={v_\varphi}(f) $$
		where
		$f=\psi^*(g) \in \mathcal{O}_X$.
		On the other hand, according to \ref{16} we have:
		\begin{align*}
			i_y \circ (\bar{\psi^*})^*(\varphi)(g)=&\left( \widetilde{(\bar{\psi^*})^*(\varphi)}\right) (g)\qquad\qquad\\
			=&\left( \varphi \circ (\overline{\psi_t^*})^*\right) \left( (g-\tilde{g})+{\mathcal J_Y}^2\right)\qquad\qquad\\ =&\varphi\left( {{{\psi}^*}}(g+{\mathcal J_Y}^2)-{{\psi}^*}(\tilde{g}+{\mathcal J_Y}^2)\right) \\
			=&\varphi\left( ({{\psi}^*}(g)-{{\psi}^*}(\tilde{g}))+{\mathcal J_Y}^2\right)\quad\quad \quad \\
			=&\varphi((f-\tilde f)+{\mathcal J_Y}^2)\\
			=&v_{\varphi}(f)
		\end{align*}
		The second equality holds
		according to \ref{5.8.0}
		and fourth equality holds
		by substituting $\tilde{\psi^*}(g)$
		for $\psi^*(\tilde{g})$.
		therefore,  up to isomorphism, the morphisms
		$d{\psi}$
		and
		$(\bar{\psi^*})^*$
		are equivalent.
	\end{proof}
	\begin{definition}
		If
		$\psi: (\tilde X, \mathcal O_X) \to (\tilde Y, \mathcal O_Y)$
		is a morphism between supermanifolds, a point
		$q\in \tilde Y$
		is said to be a regular value of
		$\psi$
		whenever at every point
		$p \in \tilde\psi^{-1}(\{q\})$,
		$\psi$ is a submersion. In other words
		$d\psi_p: T_pX \to T_qY$
		is surjective.\\
		A point	$q\in \tilde Y$ 	is said to be a critical value of
		$\psi$
		if it is not a regular value.
	\end{definition}
	
	\begin{theorem}(\textbf{Sard's theorem in supergeometry})\label{1000}
		Let $\psi: X^{m|n} \to Y^{p|q}$ be an arbitrary morphism between supermanifolds, then the set of critical values of $\psi$ i.e, $\{s\in \tilde Y|\, d\psi_t \,\,
		\text{is not }\\ \text{surjective for a t with}\,\, \tilde\psi(t)=s\}$ has (Lebesgue) measure zero.
	\end{theorem}
	\begin{proof}
		We consider the following diagram
		\begin{equation} 
			\begin{CD}
				\Xi_{X} @>(\bar{\psi^*})^*>> \Xi_{Y} \\
				@V{\pi_X}VV @VV{\pi_Y}V \\
				\tilde X @>>\tilde \psi> \tilde Y
			\end{CD}
		\end{equation}
		According to Sard's theorem in the classical case, we know that the sets of all critical values of $\bar{\psi^*}^*$, say $N$, and $\tilde\psi$ have measure zero.
		So $\pi_Y(N)$ has also measure zero in $\tilde Y$. We show that the set of all critical values of $\tilde{\psi}$ is a subset of $\pi_Y(N)\cup C_{\tilde{\psi}}$
		where $C_{\tilde{\psi}}$ is the set of critical values of $\psi$. Consider the map 
		$d\psi_t$. Set $(d\psi_t)_0= d\psi_t|_{(T_tX)_0}$
		and $(d\psi_t)_1= d\psi_t|_{(T_tX)_1}$. 
		Therefore $d\psi_t= (d\psi_t)_0 +  (d\psi_t)_1$
		and its matrix takes the form
		\begin{align}
			d \psi_{t} &= \left(\begin{array}{cc}
				(d \psi_{t})_0&0\\
				0& (d \psi_{t})_1\\
			\end{array} \right).
		\end{align}
		Now, if it is not of full rank, this means that at least one of the two blocks fails to be of full rank. Let $s$ be a critical value of $\psi$, we show that $s\in\pi_Y(N)\cup C_{\tilde{\psi}}$.
		Let $(\bar{\psi^*})^*$,
		$\Xi_X$ and $\Xi_Y$
		restricted to a coordinate neighborhood. 
		Then there is a morphism 
		\begin{align}
			(\bar{\psi^*})^*:U\times F^n &\to V \times F^q \\
			(t,v)&\mapsto (\tilde\psi(t),(\bar{\psi_t^*})^*(v))\notag.
		\end{align}
		where $t\mapsto (\bar{\psi^*})^*_t$ is a map $U\to L(F^n, F^q)$.
		Therefore, the Jacobian matrix of $\psi$ is given by
		\begin{align}
			(d (\bar{\psi^*})^*)_t  &= \left(\begin{array}{cc}
				\dfrac{\partial \tilde \psi}{\partial t}&\dfrac{\partial (\bar{\psi_t^*})^*}{\partial t}\\
				\\
				\dfrac{\partial \tilde \psi}{\partial v}=0&   \,\, \dfrac{\partial (\bar{\psi_t^*})^*}{\partial v}\\
			\end{array} \right)
		\end{align}	
		It is easily seen that the top block on the left of this matrix is equal
		to $(d\psi_t)_0$. Since  $(\bar{\psi^*})^*_t$ is a linear map of v , it is also
		seen that the lower block on  right is equal to $(\bar{\psi^*})^*_t$.
		Now diagram 5.4 implies that it is also equal to $(d\psi_t)_1$.
		If $[(d\psi_t)_0]$ fails to be 
		of full rank then
		$s\in C_{\tilde{\psi}}$.
		If $[(d\psi_t)_1]$
		fails to be of 
		full rank then 
		$s\in \pi_Y(N)$.
		This completes the proof.
	\end{proof}
	\section{Genericity}
	In the previous section, we proved that supertransversality is a stable property under small perturbations. In this section, we demonstrate a valuable fact that transversality is a generic property. This means that for any arbitrary morphism $f:X\to Y$ that is not transversal to a sub supermanifold $Z$ of $Y$, one can deform it slightly into a morphism that is supertransversal to $Z$.
	
	More precisely, we will show that there exists a supermanifold $S$ and a morphism\\$F:X \times S \to Y$ such that firstly, for some $s_0 \in \tilde S$ we have $F \circ i_{s_0} = f$, where $i_{s_0}$ is the morphism from $X$ to $X \times S$ given by the pair $(id_X, s_0)$. By $s_0$ we mean a constant morphism from $X$ to $S$. Secondly, for a given $s$ in a neighborhood of $s_0$, the morphism $F \circ i_s$ is supertransversal to $Z$. In this case, $F$ gives a smooth superhomotopy from $f$.
	\begin{theorem}({\textbf{Genericity property}})
		Let $f: X^{m|n} \to Y^{m'|n'}$ be an arbitrary morphism, and let $Z^{k|l}$ be a subsupermanifold of $Y$. There exists a smooth superhomotopy of $f$ say $f_s: X^{m|n} \to Y^{m'|n'}$ with $f_{s_0}=f$, such that for all $s$ in a sufficiently small neighborhood of $s_0$, $f_s \pitchfork Z$.
	\end{theorem}
	\begin{proof}
		We prove the theorem in four steps.
		\begin{description}
			\item[First step:] We show that for $p$ and $q$, sufficiently large, there exists an embedding
			$\phi:Y^{m'|n'} \to \mathbb{R}^{p|q}$.
			Then we consider $\psi=\phi\circ f$ a morphism from $X$ to $R^{p|q}$.
			\item [Second step:]
			We construct a smooth superhomotopy
			$\Psi:X^{m|n} \times \mathbb B^{p|q} \to \mathbb{R}^{p|q}$
			from
			$\psi$
			that is submersion. By
			$\mathbb B^{p|q}$
			we mean an open superball of
			$\mathbb R^{p|q}$.
			\item [Third step:] We first establish supertransversality
			theorem and then by using this theorem show that the morphism $\psi_s: X^{m|n} \to \mathbb R^{p|q}$
			given by
			$\psi_s=\Psi \circ i_s$. Then, by Sard's theorem in supergeometry, see theorem 5.5, for almost all
			$s \in \tilde {\mathbb B}$
			is transversal to $Z$.
			\item [Fourth step:] In the last step, it is shown that there exists a morphism
			$\pi: \mathbb R^{p|q } \to Y^{m'|n'}$ that is submersion on an open subsupermanifold of $ \mathbb R^{p|q }$ which has $Y^{m'|n'}$ as a subsupermanifold. So
			$f_s=\pi \circ \psi_s$
			is supertransversal to
			$ Z$.
		\end{description}
	\end{proof}
	\subsection{First step - Embedding theorem}
	\begin{theorem}\label{5.2.}
		Let
		$Y^{m|n}$
		be a supermanifold such that
		$\tilde{Y}$
		is compact. For
		$p|q$
		sufficiently large, there exists an embedding
		$ \phi: {Y}^{m'|n'} \to \mathbb{R}^{p|q}$.
	\end{theorem}
	\begin{proof}
		Since the
		$\tilde Y$ is compact, then we consider a finite open cover
		$ \{ \tilde U_\alpha\}_{\alpha=1}^{k}$
		for it, satisfying the following conditions: (i)
		$Y$
		is isomorphic to the superdomains in these neighborhoods. Thus, there exists isomorphisms 
		$\xi_\alpha:Y|_{U_{\alpha}} \to \bar U_{\alpha}^{m'|n'}$.
		(ii) For all
		$\alpha=1,\ldots ,k$, $\bar B_2(0) \subset\bar U_{\alpha}$
		and (iii) The family
		$\{{\tilde \xi}^{-1}_\alpha(B_1(0))\}_{\alpha=1}^k$
		is an open cover for
		$ \tilde{Y}$. If 
		$(x_i, e_j)$
		is a standard coordinate system for
		$\mathbb{R}^{m'|n'}$, then
		$\xi_\alpha^*:~\mathcal O_{\bar U_{\alpha}^{m'|n'}} \to~\mathcal O_Y(\tilde U_\alpha)$
		is as follows:
		\begin{align*}
			\xi_\alpha^*(x_i)&= f_\alpha^i \in \mathcal{O}_Y(U_\alpha)_0\\
			\xi_\alpha^*(e_j)&= g_\alpha^j \in \mathcal{O}_Y(U_\alpha)_1
		\end{align*}
		Now we consider the section
		$ \delta\in (\mathcal{O}_{R^{m'|n'}})_0$
		such that
		$\tilde \delta$
		has value 1 on
		$\bar B_1(0)$
		and vanishes outside
		$B_2(0)$
		\cite{6}.
		Then for
		$\alpha=1,\ldots ,k$
		we define sections
		$\delta_\alpha \in (\mathcal{O}_Y(\tilde U_\alpha))_0 $
		as follows:
		$$\delta_\alpha=\xi_\alpha^*(\delta) $$
		in this case
		$$\tilde \delta_\alpha(p)=\begin{cases}
			\tilde \delta(\tilde \xi_\alpha(p))\,\,\,\,\,\,\,\,\, p \in U_\alpha\\
			0\,\,\,\,\,\,\,\,\,\,\,\,\,\,\,\,\,\,\,\,\,\,\,\,\,\, p \notin  U_\alpha
		\end{cases}$$
		thus for each
		$\alpha$, set
		$\phi_\alpha: Y^{m'|n'} \to \mathbb{R}^{m'|n'}$ define as follows:
		$$\phi_\alpha^*: \mathcal{O}_{\mathbb{R}^{m'|n'}} \to \mathcal{O}_{Y^{m'|n'}}$$
		$$\phi_\alpha^* (x_i):= f_\alpha^{i} .\delta_\alpha \,, \,\,\,\,\,\, \phi_\alpha^* (e_j): = g_\alpha^{j} . \delta_\alpha$$
		therefore
		$$\tilde \phi_\alpha(p)=\begin{cases}
			\tilde \xi_\alpha(p).\tilde \delta_\alpha(p)\,\,\,\,\,\,p \in U_{\alpha}\\
			0\,\:\:\;\,\,\,\,\,\quad\qquad \,\,\, p\notin  U_{\alpha}
		\end{cases}.$$
		For convenience, we put
		\begin{align*}
			&\eta_{\alpha}^i:= \phi_\alpha^* (x_i)\\
			&\zeta_\alpha^j:= \phi_\alpha^*(e_j).
		\end{align*}
		Now we define the morphism
		$ \phi= ( \tilde{\phi} , \phi^*): Y^{m'|n'} \to \mathbb{R}^{k+km'|kn'}$
		as follows:
		\begin{align*}
			\tilde{\phi}&: \tilde{Y} \to \mathbb{R}^{k+km}\\
			&\,\,\,\,\,\,\, p \mapsto (\tilde\delta_1(p),\cdots,\tilde\delta_k(p), \tilde \phi_1(p) ,\cdots, \tilde \phi_k(p))
		\end{align*}
		and  according to the coordinates
		$(x_i,e_j)$
		for
		$\mathbb{R}^{k+km'|kn'}$,
		$\phi^*: \mathcal{O}_{\mathbb{R}^{k+km'|kn'}} \to \mathcal{O}_{Y^{m'|n'}}$
		as follows:\\			
		\begin{tabular}{c c}
			$\begin{cases}
				x_i \mapsto \delta_i=\eta_0^i\,\,\,\,\,\,\,\, 1 \le i \le k\\
				x_i \mapsto \eta_1^i \,\,\,\,\,\,\,\, k+1 \le i \le k+m'\\
				x_i \mapsto \eta_2^i \,\,\,\,\,\,\,\, k+m'+1 \le i \le k+2m'\\
				\vdots\\
				x_i \mapsto \eta_k^i \,\,\,\,\,\,\,\, k+(k-1) m'+1 \le i \le k+km'
			\end{cases}$
			& 
			$\begin{cases}
				e_j\mapsto \zeta_1^j \,\,\,\,\,\,\,\,\,\, 1 \le j \le n\\
				e_j \mapsto \zeta_2^j \,\,\,\,\,\,\,\,\,\,n'+1 \le j \le 2n\\
				\vdots\\
				e_j \mapsto \zeta_k^j \,\,\,\,\,\,\,\,\,\, (k-1)n'\le j \le kn'\\
			\end{cases}$
		\end{tabular}\\
		we show that the rank of
		$\phi$
		at point $p$ is equal to
		$m|n$. For this end, let $p$ be a point in $\tilde \xi_{\alpha_0}^{-1}(B_1(0))$ and let $(y_a, \theta_b)$ and $(x_i, e_j)$ be coordinate systems on $U_{\alpha}$ and $R^{k+ km'|kn'}$ respectively. Then one has
		$$[ J \phi ]=  \begin{pmatrix}
			\dfrac{\partial \eta_\alpha^i }{\partial y_a}&-\dfrac{\partial \eta_\alpha^i}{\partial \theta_b}\\
			\\
			\dfrac{\partial \zeta_\alpha^j }{\partial y_a}& \,\, \dfrac{\partial \zeta_\alpha^j}{\partial \theta_b}\\
		\end{pmatrix}
		$$
		where 
		$1 \le a \le m'$
		,
		$1 \le b \le n'$
		,
		$1 \le i \le k+km'$
		and
		$1 \le j \le kn'$
		.
		The matrix
		$[J\phi]$
		induces the following matrix:
		\begin{equation}\label{5.1}
			\tilde{[J\phi]}_{(p)}= 
			\begin{pmatrix}
				\dfrac{\partial \eta_\alpha^i }{\partial y_a}^\sim(p)&0\\
				0 & \dfrac{\partial \zeta_\alpha^j}{\partial \theta_b}^\sim(p)\\
			\end{pmatrix}.
		\end{equation}
		Since
		\begin{align*}
			&\eta_{\alpha}^i= f_\alpha^{i} .  \delta_\alpha\Rightarrow\partial\eta_{\alpha}^i=\partial f_\alpha^{i} .  \delta_\alpha+f_\alpha^{i}\partial \delta_\alpha\\
			&\zeta_\alpha^j= g_\alpha^{j} . \delta_\alpha\Rightarrow\partial\zeta_{\alpha}^i=\partial g_\alpha^{j} .  \delta_\alpha+(-1)^{p(g_\alpha^i)}g_\alpha^{j}\partial \delta_\alpha
		\end{align*}
		then it can be seen that the matrix \ref{5.1} has a submatrix as follows:
		$$ \begin{pmatrix}
			\dfrac{\partial f^i_{\alpha_0} }{\partial y_a}(p)&0\\
			0 & \dfrac{\partial g_{\alpha_0}^j}{\partial \theta_b}(p)\\
		\end{pmatrix} $$
		where
		$1 \le i,a \le m$
		and
		$ 1\le j,b \le n$
		thus its rank is equal to
		$m'|n'$. So, the matrix
		\ref{5.1}
		has at least a rank
		$m'|n'$. Since the number of columns in the block
		$\left( \dfrac{\partial \eta_\alpha^i }{\partial y_a}\right)$
		at point
		$p$
		is equal to
		$m'$ and the number of columns in the block
		$ \left( \dfrac{\partial \zeta_\alpha^j}{\partial \theta_b}\right)$
		at the point
		$p$
		is equal to
		$n'$
		, the rank of the matrix
		$(J\phi)$
		cannot be greater than
		$m'|n'$ so
		$$ rank [J \phi (p) ]= m'|n'.$$
		Since
		$\tilde{\phi}: \tilde{Y}^{m'}\to \mathbb{R}^{k+km'}$ is an embedding therefore
		$\phi$
		is an embedding too.
	\end{proof}	
	For convenience set $\psi:=\phi\circ f$.
	\subsection{Second step - Constructing a smooth homotopy}
	\begin{definition}
		\label{def 1-3-0}
		Let
		$\tilde{\mathbb B}$
		be an open ball in
		$ \mathbb{R}^{p}$
		then we define an open superball as follows:
		$$\mathbb B^{p|q}:=(\tilde {\mathbb B}, \mathcal{O}_{\mathbb{R}^{p|q}}(\tilde {\mathbb B})) $$.
	\end{definition}
	\begin{definition}
		Suppose
		$\psi: X^{m|n}\to \mathbb{R}^{p|q}$
		be a morphism between supermanifolds. If
		$\mathbb B^{p|q}$
		be an open superball then, we define a smooth superhomotopy 
		$\Psi=(\tilde\Psi,\Psi^*): X^{m|n} \times \mathbb B^{p|q} \to \mathbb{R}^{p|q}$
		as follows:
		\begin{equation*}
			\begin{array}{lccccc}
				&\tilde\Psi: \tilde X\times\tilde{\mathbb B} \to \mathbb{R}^p \\
				&\,\,\,\,\,\,\,\,\,\,\,\,\,\,\,\,\qquad(x,s) \mapsto \tilde \psi(x)+s\\
				\\
				&\Psi^*:\mathcal{O}_{\mathbb{R}^{p|q}} \to \mathcal{O}_{X^{m|n}\times \mathbb B^{p|q}}\\
				&\,\,\,\,\,\,\,\,\,\,\,\,\,\,\,\,\,\,\,\, x_i \mapsto \psi^*(x_i)+x_i \\
				&\,\,\,\,\,\,\,\,\,\,\,\,\,\,\,\,\,\,\,\, e_j \mapsto  \psi^*(e_j)+e_j
			\end{array}
		\end{equation*}
		where $(x_1, ..., x_p; e_1, ..., e_q)$ is a global coordinate on $\mathbb R^{p|q}$ and $\psi^*$ is the pullback $\mathcal{O}_{\mathbb R^{p|q}}\to \mathcal{O}_{X^{m|n}}$.
	\end{definition}
	\begin{proposition}
		$\Psi:X^{m|n} \times \mathbb B^{p|q} \to \mathbb{R}^{p|q}$ is a submersion.
	\end{proposition}
	\begin{proof}
		Let $(\xi_1, ..., \xi_m; \eta_1, ..., \eta_n)$ be a coordinate system around the arbitrary  point $x_0$ in $\tilde{X}$. We need to examine the Jacobian matrix of morphism $\Psi$.
		For convenience, we show even coordinates
		$(\xi_1,...,\xi_m,x_1,...,x_p)$
		with
		$y_l$
		and odd coordinates
		$(\eta_1,...,\eta_n,e_1,...,e_q)$
		with
		$\theta_k$. Therefore, $[J\Psi]$ decomposed into four blocks as follows:
		$$[ J \Psi ]=  \begin{pmatrix}
			\dfrac{\partial ( \psi^*(x_i)+x_i) }{\partial y_l}& - \dfrac{\partial  ( \psi^*(x_i)+x_i)}{\partial \theta_k}\\
			\\
			\dfrac{\partial  ( \psi^*(e_j)+e_j) }{\partial y_l}&  \dfrac{\partial  (\psi^*(e_j)+e_j)}{\partial \theta_k}\\
		\end{pmatrix}.$$
		This matrix induces the following matrix:
		$$[ J \tilde\Psi(x_0,s_0) ]= \begin{pmatrix}
			\left(\dfrac{\partial ( \psi^*(x_i)+x_i) }{\partial y_l}\right)^\sim(x_0,s_0)&0\\
			0 & \left(\dfrac{\partial( \psi^*(e_j)+e_j)}{\partial \theta_k}\right)^\sim(x_0,s_0)\\
		\end{pmatrix}$$
		Where, regarding the definition of $\Psi$, one has the following equalities for non-zero blocks are as follows:
		$$
		\left(\dfrac{\partial ( \psi^*(x_i)+x_i) }{\partial y_l}\right)^\sim(x_0,s_0)=(*\,| I_p)$$
		and
		$$ \left(\dfrac{\partial( \psi^*(e_j)+e_j)}{\partial \theta_k}\right)^\sim(x_0,s_0)= (*\,| I_q),$$
		where
		$I_p$
		and
		$I_q$
		are identity matrices
		$p\times p$
		and
		$q\times q$, respectively.
		Thus
		$\Psi$
		is a submersion.
	\end{proof}
	\subsection{Third step - Supertransversality Theorem}
	\begin{lemma} \label{6.6.0}
		Let $s \in \tilde{\mathbb B}$	and 
		$i_s : X \to X \times \mathbb B$
		is the morphism given by a pair
		$(id_X, s)$, then for $w \in T_xX$ we show that
		$di_s(w)=(w,0)$.
	\end{lemma} 
	\begin{proof}
		Let
		$ i_s=( \tilde{i}_s, i_s^*)$
		where
		$ \tilde{i}_s: \tilde{X} \to \tilde{X} \times \tilde{\mathbb B}^{p|q}$. By considering the coordinates\\
		$(\xi_1,...,\xi_m,\zeta_1,...,\zeta_n)$ for neighborhood 
		$U$ around $x$ on
		$X$
		and coordinates
		$(x_1,...,x_p, e_1,...,e_q)$
		for neighborhood 
		$V$ around $s$ on
		$ \mathbb B$,
		the morphism $i_s^*: \mathcal{O}_{X \times \mathbb B^{p|q} }(\tilde U \times \tilde V) \to \mathcal{O}_X(\tilde U)$ can be written in these coordinates as follows:
		\begin{align*}
			\xi_i \mapsto \xi_i \,\,\,\,\,\,\,\,\,\,\,\, 1 \le i \le m\\
			\zeta_j \mapsto \zeta_j\,\,\,\,\,\,\,\,\,\,\, 1 \le j \le n\\
			x_a \mapsto \tilde{x}_a(s) \,\,\, 1 \le a \le p\\
			e_b \mapsto 0 \,\,\,\,\,\,\,\,\,\,\,\,\,\, 1 \le b \le q
		\end{align*}
		Since
		$\left( \dfrac{\partial }{\partial \xi_i}\right) _x , \left( \dfrac{\partial }{\partial \zeta_j}\right) _x$ is a basis for $T_xX$
		thus, each
		$w \in T_xX$
		can be written as
		$w= w^i \partial / \partial \xi_i +\theta^j \partial / \partial \zeta_j $. Also $T_{(x,s)} ( X \times \mathbb B^{p|q})$ has a basis as
		$ \left( \dfrac{\partial }{\partial \xi_i}\right) _{(x,s)} , \left( \dfrac{\partial }{\partial x_a}\right)_{(x,s)}, \left( \dfrac{\partial }{\partial \zeta_j}\right)_{(x,s)},\left( \dfrac{\partial }{\partial e_k}\right)_{(x,s)} $ therefore, we have
		\begin{align*}
			& di_s(w)= k_{1i}\dfrac{\partial }{\partial \xi_i} +k_{2a} \dfrac{\partial }{\partial x_a}+ k_{3j} \dfrac{\partial }{\partial \zeta_j}+  k_{4k}\dfrac{\partial }{\partial e_k} \\
			d i_s (w) (\xi_{i_0}) & =k_{1{i_0}} \\
			&= w( i_s^* (\xi_{i_0}))= w(\xi_{i_0})= ( w_i \partial / \partial \xi_i + \theta_j \partial / \partial \zeta_j )(\xi_{i_0})=w_{i_0}\\
			di_s(w) (\zeta_{j_0}) & = k_{3{j_0}}\\
			&= w( i_s^* (\zeta_{j_0}))= w(\zeta_{j_0})= ( w_i \partial / \partial \xi_i + \theta_j \partial / \partial \zeta_j )(\zeta_{j_0})=\theta_{j_0}\\
			di_s(w) (x_a) &= k_{2a} = w( i_s^*(x_a))= w( \tilde{x}_a(s))=0\\
			di_s(w) (e_b) &= k_{4k} = w( i_s^*(e_b))= w( 0)=0.
		\end{align*}
		Thus $di_s(w)=(w,0)$.
	\end{proof}
	\begin{proposition}\label{p5.7}
		Suppose
		$\Psi: X^{m|n} \times \mathbb B^{p|q} \to \mathbb{R}^{p|q}$
		be a morphism between supermanifolds, and
		$Z^{k|l}$ be a subsupermanifold of
		$ \mathbb{R}^{p|q}$ such that
		$ \Psi\pitchfork Z$.
		Then for almost all
		$s \in \tilde{\mathbb B}$,
		$ \psi_s \pitchfork Z$
		where
		$\psi_s = \Psi \circ i_s$
		and
		$i_s : X \to X \times \mathbb B$
		is the morphism as in the lemma
		\ref{6.6.0}.
	\end{proposition}
	\begin{proof}
		By using theorems \ref{3.5.} and \ref{1000}, the proof is similar to
		the proof of transversality theorem \cite{3}.
	\end{proof}
	\subsection{Fourth step - A morphism from \texorpdfstring{$\mathbb R^{p|q}$}- to an embedded subsupermanifold}
	In this section, we introduce a morphism that is needed to prove the genericity property of supertransversality.	
	\begin{proposition}
		There exists a morphism
		$\pi:\mathbb R^{p|q} \to Y^{m'|n'} $
		which is submersion at each point of an open neighborhood of $\tilde{Y}$ in $\mathbb R^p$.
	\end{proposition}
	\begin{proof}
		By theorem \ref{5.2.} there exists an embedding $\phi:Y^{m'|n'} \to \mathbb R^{p|q}$. Thus, for each point
		$y \in \tilde Y$
		there exist open neighborhoods
		$\tilde U_\alpha$
		in
		$\tilde Y$
		and
		$\tilde V_\alpha$
		in
		$\mathbb R^p$ and a coordinate system $(r_1, ..., r_p; e_1, ..., e_q)$ such that the following homomorphism is surjective:
		$$\phi^*_{\tilde V_{\alpha}}: \mathcal{O}_{\mathbb R^{p|q}}(\tilde V_{\alpha})\to \mathcal{O}_{Y^{m'|n'}}(\tilde U_{\alpha})$$
		and $\phi^*_{\tilde V_{\alpha}}(g)= f$ where 
		$$f(r_1, ..., r_{m'}; e_1, ..., e_{n'}):=g(r_1, ..., r_{m'}, 0, ..., 0; e_1, ..., e_{n'}, 0, ..., 0).$$
		Now, we consider the open covering
		$\left\lbrace \tilde V_\alpha\right\rbrace _\alpha \bigcup\left\lbrace \mathbb R^p/\tilde Y\right\rbrace $
		for
		$\mathbb{R}^{p}$ and choose a partition of unity
		$\{\zeta_\alpha\}_\alpha $
		corresponding to this covering where for
		$1\leq \alpha\leq k$,
		$supp\,\zeta_\alpha \subset \tilde V_\alpha$
		and
		$supp\,\zeta_0 \subset \mathbb R^p/\tilde Y $\cite{4,5}. For any arbitrary open neighborhood
		$\tilde W$
		in
		$\tilde Y$
		define
		\begin{equation}\label{5.18.}
			\begin{cases}
				\pi^*:&\mathcal{O}_{Y^{m'|n'}}(\tilde W)\to \mathcal{O}_{\mathbb R^{p|q}}(\tilde \pi^{-1}(\tilde W))\\ 
				& \,\,\,\,\,\,\,\,\,\,\,\,\quad\, f \mapsto \sum_\alpha g_{\alpha}.\zeta_{\alpha}
			\end{cases}
		\end{equation}
		where $g_{\alpha}$ is an element of $\mathcal{O}_{\mathbb R^{p|q}}(\tilde V_{\alpha})$, for a proper open subset of $\tilde V_\alpha$, such that $$\phi^*_{\tilde V_{\alpha}}(g_{\alpha})=f|_{\tilde W\cap \tilde U_{\alpha}}.$$
		Let $\tilde V=\cup \tilde V_{\alpha}$. We show that
		$$d\pi_y: T_y\mathbb R^{p|q} \to T_{\tilde \pi(y)}Y^{m'|n'}$$
		is surjective for each $y\in \tilde V$.\\
		For each $\tilde \pi(y)\in \tilde Y$ there exists an index $\alpha$ such that $\tilde \pi(y) \in \tilde U_\alpha$ and $y \in \tilde V_\alpha$ therefore according to above there exists coordinate system $(r_1, ..., r_p; e_1, ..., e_q)$ for $ V_\alpha$ and $(r_1,...,r_{m'};e_1,...,e_{n'})$ for $U_\alpha$ such that the tangent space $T_y{\mathbb R^{p|q}}$ has $\left( \dfrac{\partial}{\partial r_i}\right)_y,\left( \dfrac{\partial}{\partial e_j}\right)_y$ as a basis. 
		With respect to \ref{5.18.} for $1 \leq i \leq m'$ we have
		$$\pi_y^*(r_{i})=\sum_\alpha r_{i}.\zeta_a=r_i$$
		and similarly  
		$\pi_y^*(e_j)=e_j$ for $1 \leq j \leq n'$.
		So the Jacobian matrix of morphism $\pi$ is as follows:
		$$
		[J\pi]_y=
		\begin{blockarray}{ccc}
			p & q  \\
			\begin{block}{(cc)c}
				I_{m'} & 0 & m' \\
				0 & I_{n'} & n' \\
			\end{block}
		\end{blockarray}
		$$
		Thus $[J\pi]_y$ is a full rank matrix and $d\pi_y$ is surjective.
	\end{proof}
	\begin{theorem}
		$f_s=\pi \circ \psi_s:X^{m|n} \to Y^{m'|n'} $ is
		supertransversal to
		$Z^{k|l}$.
	\end{theorem}
	\begin{proof}
		By proposition \ref{p5.7} for one $s\in \tilde {\mathbb B}^p$, $\psi_s\pitchfork Z$. 
		Thus, for any $y \in \tilde Z$ we have
		$$Im\,d\psi_s+T_yZ=T_y\mathbb R^{p|q}$$
		where $\tilde \psi_s(x)=y$. So, for any arbitrary tangent vector $a \in T_y\mathbb R^{p|q} $ there exists $u \in T_xX$  such that 
		\begin{equation}\label{5.19}
			d\psi_s(u)-a \in T_yZ.
		\end{equation}
		According to the previous proposition for the neighborhood $\tilde V$ in $\mathbb R^{p}$
		$$d\pi_y: T_yV \to T_yY$$
		is surjective, thus for any vector $b\in T_yY$ one has $b=d\pi(a)$ for some $a \in T_yV$.
		By applying $d\pi_y$ to \ref{5.19} we have
		$$d\pi_y(d\psi_s(u)-a)\in T_yZ$$
		then 
		$$d\pi_y d\psi_s(u)-b \in T_yZ$$
		where $u \in T_x(\psi^{-1}_s(Z))$, this means $f_s \pitchfork Z$.
	\end{proof}											\section{ \texorpdfstring{$\Pi$-}-symmetric supermanifolds} \label{5.4.0}
	In this section, we will introduce the category of 
	$\Pi $-symmetric supermanifolds. and then examine supertransversality in this category. At first, we talk about the concepts of free modules and
	$\Pi$-symmetries \cite{2, 11, 12}.
	A $\Pi$-symmetric free
	$A$-module $S$ of rank $m|m$, is defined to be an $A$-module isomorphic to
	$$ A^m \oplus (\Pi A)^m.$$
	Now we consider
	$A$-module isomorphism
	$ P_S: S \to \Pi S$ such that
	$P_S^2=id_S$
	where
	$\Pi$
	be a parity change functor.
	We call
	$P_S$, a
	$\Pi$ -symmetry on
	$S$.
	Objects in the category of $\Pi$-symmetric supermanifolds are supermanifolds such that there exists a $\Pi$-symmetry on their tangent sheaf. From now, for convenience, by $\Pi$-symmetry on a supermanifold we mean a $\Pi$- symmetry on its tangent sheaf.
	\subsection{Category of $\Pi$-symmetric supermanifolds}
	Let $M$ and $N$ be $\Pi$- symmetric supermnifolds. Then the set of morphisms $Hom(M, N)$  consists of all morphism $\psi: M\to N$ such that the following diagram commutes
	\begin{equation}\label{5.4.}
		\begin{CD}
			\psi_*(\mathcal T_M)@> {\bar D(\psi)}>>\mathcal T_{N,M}@.\qquad \\
			@ V{P_M}VV @ VV{\bar P_N} V @.\\
			\Pi\psi_*(\mathcal T_M)@>>{\bar D(\psi)^\Pi}>\Pi\mathcal T_{N,M}@.
		\end{CD}
	\end{equation}
	where $P_M$ is $\Pi$-symmetry on $\mathcal{T}_M$ and $\bar{P}_N$ is a  morphism induced by $\Pi$-symmetry $P_N$ on $\mathcal{T}_{N,M}:= Der(\mathcal{O}_N \to \mathcal{O}_M)$ such that $\bar{P}_N(f\partial/\partial y):=fP_N(\partial/\partial y)$ for each $f \in \mathcal{O}_M$ and 
	$\partial/\partial y$ is a coordinate vetor field on $N$.
	
	In addition $\bar{D}(\psi): \psi_*(\mathcal{T}_M)\to \mathcal{T}_{N,M}$ is defined as follows:
	$\bar{D}(\psi){X}(g)= X(\psi^*(g))$ where $X \in Der(\mathcal O_M)$ and $g \in \mathcal O_N$. The next proposition is a straightforward consequence of the above definitions.
	\begin{proposition}
		$\Pi$- symmetric supermanifolds and the morphisms between them constitute a category denoted by $\mathscr{S}_\Pi$.
	\end{proposition}
	For local consideration, if
	$(x_1,... ,x_m)$
	even coordinates on
	$M$
	around the point
	$q$, then
	$(T_qM)_0$
	is freely generated by
	$\left\lbrace\dfrac{\partial}{\partial x_1},...,\dfrac{\partial}{\partial x_m}\right\rbrace$
	so
	$\left\lbrace P_M\left(\dfrac{\partial}{\partial x_1}\right),...,P_M\left(\dfrac{\partial}{\partial x_m}\right)\right\rbrace $
	is a basis for
	$(T_qM)_1.$ 
	\begin{proposition}
		Suppose
		$S=S_0 \oplus S_1$ be a $\Pi$-symmetric free module.
		For any even module's homomorphism  
		$P_S: S \to \Pi S$, provided that
		$P_S^2=id_S$,
		there exists an odd  module's homomorphism
		$P_\Pi:S \to S$
		such that
		$P_\Pi^2=id_S$
		and vice versa.
	\end{proposition}
	\begin{proof}
		First, note that
		\begin{align*}
			& ( \Pi S)_0= \{ \Pi a, a \in S_1 \}\\
			& ( \Pi S)_1= \{ \Pi a, a \in S_0 \} .
		\end{align*}
		Now set 
		$ \Pi: \Pi S \to S$,
		we define
		$$P_\Pi:= \Pi \circ P_S : S \to S$$
		Since
		$$P_S \circ \Pi= \Pi \circ P_S $$
		one has
		\begin{align}		P_{\Pi}^2=(\Pi \circ P_S) \circ ( \Pi \circ P_S) &= \Pi \circ \Pi \circ P_S \circ P_S\\
			&= \Pi^2 \circ P_S^2= id_S.
		\end{align}
		Since
		$P_S: S \to \Pi S$
		is an isomorphism, therefore put for even coordinates
		$e_i$
		in
		$S_0$, one has
		$$ P_S e_i \in ( \Pi S)_0= \{ \Pi a| a \in S_1\}$$
		thus
		$S=\left\langle e_1,...,e_m, P_\Pi e_1,...,P_\Pi e_m\right\rangle $.
		For convenience, $P_{\Pi}$ also called $\Pi$-symmetry on $S$.
	\end{proof}
	
	\begin{proposition}\label{9.5.}					
		Suppose
		$S$
		be a $\Pi$-symmetric free module with the basis
		$$S =< e_1,...,e_m, P_\Pi e_1,...,P_\Pi e_m>$$
		where $P_{\Pi}$ is an $\Pi$-symmetry on $S$. Then
		$A \subset S$
		is a
		$\Pi$-symmetric submodule if and only if for each element
		$v= x^{i}e_i+\xi^iP_\Pi e_i$ in $A$, the element
		$v_\Pi=-\xi^ie_i+ x^{i}P_\Pi e_i$
		is also in $A$.
	\end{proposition}
	\begin{proof}
		According to the previous proposition, we have
		\begin{align*}
			&x^i e_i + \xi^i   P_\Pi e_i \in A\\
			& \Rightarrow P_S(x^i e_i+ \xi^i  P_\Pi e_i)  \in \Pi A \\
			& =x^i P_S e_i + \xi^i P_S  P_\Pi e_i  \in \Pi A\\
			& =x^i \Pi  P_\Pi e_i + \xi^i \Pi P_\Pi  P_\Pi e_i  \in \Pi A\\
			& =x^i \Pi  P_\Pi e_i + \xi^i \Pi e_i \in \Pi A\\
			& =\Pi ( x^i  P_\Pi e_i) + \Pi ( -\xi^i) e_i \in \Pi A\\
			& =\Pi ( x^i P_\Pi e_i - \xi^i e_i) \in \Pi A\\
			& \Rightarrow x^i  P_\Pi e_i - \xi^i e_i \in A
		\end{align*}
	\end{proof}
	\subsection{$\Pi$-symmetricity of preimage}
	In the following, we work in the category of $\Pi$-symmetric supermanifolds. 
	\begin{proposition}\label{5.12.}
		Open superball 
		$ \mathbb B^{n|n} =( \tilde{\mathbb B}, \mathcal{O}_{\mathbb{R}^{n|n}} ( \tilde{\mathbb B}))$
		in
		$ \mathbb{R}^{n|n}$,
		is
		$ \Pi$-symmetric.
	\end{proposition}
	According to the previous proposition, the proof is obvious.
	\begin{proposition}
		If
		$X$
		and
		$\mathbb B$
		be $\Pi$-symmetric supermanifolds then
		$X \times \mathbb B$ is an object of $\mathscr{S}_\pi$.
	\end{proposition}
	\begin{proof}
		According to the assumption, we consider odd homomorphisms
		$P_X:\mathcal T_X \to \mathcal T_X$
		and
		$P_\mathbb B: \mathcal T_\mathbb B \to \mathcal T_\mathbb B$
		with conditions
		$P_X^2=id_{\mathcal T_X}$
		and
		$P_\mathbb B^2=id_{\mathcal T_\mathbb B}.$
		Since the
		$ \mathcal T _{X\times \mathbb B}=\mathcal T_X \times \mathcal T_\mathbb B$
		( In fact, we have
		$ \mathcal T_{X\times \mathbb B}= \mathcal T_X \oplus \mathcal T_\mathbb B $),
		Put
		$P_{X\times \mathbb B}:=P_X \times P_\mathbb B.$
		We have
		$$P_{X\times \mathbb B}^2=P_X^2 \times P_\mathbb B^2=id_{\mathcal{T}_X}\times id_{\mathcal{T}_\mathbb B} =id_{\mathcal{T}_ {X\times \mathbb B}}.$$
		therefore
		$X\times \mathbb B$
		is $ \Pi$-symmetric.	
	\end{proof}

	Now, suppose
	$\Psi: X^{m|m} \times \mathbb B^{n|n} \to \mathbb R^{n|n}$
	is supertransversal to
	$Z^{k|k}$. By theorem \ref{3.5.}, $W=\Psi^{-1}(Z)$
	is subsupermanifold of $X\times B^{n|n}$. We show that
	$W$ is
	$\Pi$-symmetric. For this purpose, it is necessary to consider a local form of the morphisms
	$\bar{D}(G_U \circ \Psi)$ and $\bar{D}(i_W)$, c.f. \ref{5.4.} for their definitions, as follows:
	$$\bar{D}( G_U \circ \Psi): Der(\mathcal{O}_{X \times \mathbb B}) \to Der(\mathcal{O}_{ \mathbb{R}^{n-k| n-k}}) \otimes \mathcal{O}_{X \times \mathbb B}$$
	$$\bar D(i_W): Der(\mathcal O_{W^{m+k|m+k}}) \to Der(\mathcal O_{X^{m|m} \times \mathbb B^{n|n}}) \otimes \mathcal O_W.$$
	According to the theorem
	\ref{2.10.} there exists a local coordinate system
	$(r_{1},...,r_{2n-2k})$
	around
	$0 \in~{\mathbb R^{n-k}}$
	and a coordinate system
	$(e_1,...,e_{2m+2k}, r_1 ,..., r_{2n-2k})$
	for
	$ X^{m|m}\times \mathbb B^{n|n}$
	at a point of
	$\widetilde{(G_U \circ\Psi)}^{-1}(0)$
	such that
	\begin{equation} \label{6}
		\begin{cases}
			G_U \circ\Psi: X^{m|m} \times {\mathbb B}^{n|n} \to \mathbb{R}^{n-k|n-k}\\
			(G_U \circ\Psi)^*:\mathcal{O}_{\mathbb{R}^{n-k|n-k}} \to \mathcal{O}_{ X^{m|m} \times \mathbb B^{n|n}}\\
			\;\;\;\;\qquad\qquad\qquad\quad\;	r_i \mapsto r_i \;\;\;\;\; 1\;\leq i  \leq\;2n-2k
		\end{cases}.
	\end{equation}
	Thus
	$(e_1,...,e_{2m+2k})$
	is a coordinate system for
	$ \mathcal{O}_{W}$
	and pullback map of the morphism
	$i_W: W \to X \times \mathbb B^{n|n}$
	is as follows:
	\begin{equation} \label{7}
		\begin{cases}
			i_W^*:\mathcal{O}_{ X^{m|m} \times \mathbb B^{n|n}} \to \mathcal{O}_{W^{m+k|m+k}} \\
			\;\;\;\;\qquad\qquad\quad	e_i \mapsto e_i \;\;\;\;\; 1\;\leq i  \leq\;2m+2k\\
			\;\;\;\;\qquad\qquad\quad	r_j \mapsto 0 \;\;\;\;\;\; 1\;\leq j  \leq\;2n-2k
		\end{cases}.
	\end{equation}
	Therefore, vector fields
	$\dfrac{\partial}{\partial e_i},\dfrac{\partial}{\partial r_j} $
	for
	$ 1 \le i \le 2m+2k$
	and
	$ 1 \le j \le 2n-2k$
	are local generators for module
	$ Der ( \mathcal{O}_{X \times \mathbb B^{n|n}}),$
	and
	$\dfrac{\partial}{\partial e_i}$
	for
	$ 1 \le i \le 2m+2k$
	are also local generators for module
	$ Der ( \mathcal{O}_W)$.
	So for a
	$Y \in Der ( \mathcal{O}_W)(\Psi^{-1}(U) \cap W)$
	we have
	$$Y=\sum_{i=1}^{2m+2k}f_i\dfrac{\partial}{\partial e_i} \qquad f_i \in \mathcal{O}_W. $$
	Thus
	$\bar{D}(i_W)Y$
	is as follows:
	\begin{equation}\label{5.10}
		\bar{D}(i_W)Y=\sum_{i=1}^{2m+2k}\dfrac{\partial}{\partial e_i} \otimes_{\mathcal{O}_{X \times \mathbb B^{n|n}}}^{i_W^*} d e_i(Y)
	\end{equation}
	where
	$ d e_i \in \Omega^1(\mathcal{O}_W)$
	and
	$d e_i(\dfrac{\partial}{\partial e_j})=\delta_{ij}$
	so
	$d e_i(Y)=f_i$. It is worth noting that $\bar{D}(i_W)Y$ belong to   $Der(\mathcal O_{X\times \mathbb B^{n|n}}, \mathcal{O}_W)$ and this sheaf, locally up to isomorphism, is equal to \\ $Der(\mathcal O_{X\times \mathbb B^{n|n}})\otimes_{\mathcal{O}_{X\times \mathbb B^{n|n}}}~{\mathcal{O}_W}$.
	In addition, for each
	$T \in Der ( \mathcal{O}_{X \times \mathbb B^{n|n}})(U)$
	we have
	$$T=\sum_{j=1}^{2m+2k}g_j\dfrac{\partial}{\partial e_j}+\sum_{l=1}^{2n-2k}h_l\dfrac{\partial} {\partial r_l}.$$
	So
	\begin{equation}\label{5.11}
		\bar{D}(G_U \circ \Psi)T = \sum_{l=1}^{2n-2k}\dfrac{\partial}{\partial r_l} \otimes_{\mathcal{O}_ {\mathbb R^{n-k|n-k}}}^{(G_U \circ \Psi)^*} h_l.
	\end{equation}
	It was also shown that
	$G_U\circ \Psi$
	is submersion	
	\begin{proposition}
		According to the equalities (\ref{5.10}) and (\ref{5.11}), we have
		$$\bar{D}(i_W)(\mathcal{T}_W(\Psi^{-1}(U)))= ker\, (\bar{D}( G_U \circ \Psi)\otimes 1_ {\mathcal{O}_W})$$
		where
		$U$
		is an open neighborhood  in
		$\mathbb R^{n|n}$
		and
		$G_U: U\to \mathbb R^{n-k|n-k}$
		is a morphism such that
		$Z\cap U=G_U^{-1}(0)$
		and
		$W=\cup \Psi^{-1}(Z \cap U)$.
	\end{proposition}
	\begin{proof}
		We start from the left side of the equality.\\ According to the morphism
		$i_W: W \to~X \times~\mathbb B^{n|n}$
		we have
		$$\bar{D}(i_W): Der (\mathcal{O}_W) \to Der(\mathcal{O}_{X \times \mathbb B^{n|n}})\otimes \mathcal { O}_W$$
		where
		$\mathcal{T}_W(\Psi^{-1}(U))=Der (\mathcal{O}_W)(\Psi^{-1}(U)).$
		Consider the following composition
		$$ Der(\mathcal{O}_{W} )\xrightarrow{\bar{D}(i_W)} Der(\mathcal{O}_{ X \times \mathbb B^{n|n} }) \otimes_{\mathcal{O}_{X \times \mathbb B^{n|n}}}^{i_W^*} \mathcal{O}_W \xrightarrow{\bar{D}(G_{U}\circ \Psi) \otimes1_{\mathcal{O}_{W}}} Der( \mathcal{O}_{ \mathbb{R}^{n-k|n-k} })\otimes_{\mathcal{O}_{\mathbb R^{n-k|n-k}}}^{(G_U \circ \Psi)^*} \mathcal{O}_{X \times \mathbb B^{n|n}}\otimes_{\mathcal{O}_{X \times \mathbb B^{n|n}}}^{i_W^*} \mathcal{O}_{W}$$
		we will show for each
		$Y \in Der(\mathcal{O}_{W})(\Psi^{-1}(U))$,
		$(\bar{D}(G_{U}\circ \Psi) \otimes1_{\mathcal{O}_{W}})(\bar{D}(i_W)Y)=0$.
		By \ref{5.10} one has
		\begin{align*}	
			&(\bar{D}(G_{U}\circ \Psi) \otimes1_{\mathcal{O}_{W}})(\bar{D}(i_W)Y)=\\
			&(\bar{D}(G_{U}\circ \Psi) \otimes1_{\mathcal{O}_{W}})(\sum_{i=1}^{2m+2k}\dfrac{\partial}{\partial e_i} \otimes_{\mathcal{O}_{X \times \mathbb B^{n|n}}}^{i_W^*} f_i)=\\
			&\sum_{i=1}^{2m+2k}\bar{D}(G_U \circ \Psi)\dfrac{\partial}{\partial e_i} \otimes_{\mathcal{O}_{X \times \mathbb B^{n|n}}}^{i_W^*}\;1_{\mathcal{O}_W}f_i=\\
			&\sum_{i=1}^{2m+2k}\left(\sum_{l=1}^{2n-2k}\dfrac{\partial}{\partial r_l} \otimes_{\mathcal{O}_{\mathbb R^{n-k|n-k}}}^{(G_U \circ \Psi)^*} 0\right)   \otimes_{\mathcal{O}_{X \times S}}^{i_W^*}f_i =0
		\end{align*}
		as a result,
		$$\bar{D}(i_W)(\mathcal{T}_W(\Psi^{-1}(U))
		\subseteq
		ker\, (\bar{D}( G_U \circ \Psi)\otimes 1_{\mathcal{O}_W}).$$	
		On the other hand, for the morphism
		$\bar{D}(G_U \circ \Psi)\otimes 1_{\mathcal{O}_W}$
		we have
		$$dim (\, Der (\mathcal{O}_{X \times \mathbb B^{n|n}})\otimes \mathcal{O}_W)=dim \,ker (\bar{D}( G_U \circ \Psi)\otimes 1_{\mathcal{O}_W})+ dim\, Im \,(\bar{D}(G_U \circ \Psi)\otimes 1_{\mathcal{O}_W}).$$
		Since
		$dim( Der (\mathcal{O}_{X \times \mathbb B^{n|n}})\otimes\mathcal{O}_W)= (2m+2n)$ and $G_U\circ \Psi$ is submersion, one has
		$$ dim\, Im(\bar{D}(G_U \circ \Psi)\otimes 1_{\mathcal{O}_W})=dim (Der( \mathcal{O}_{ \mathbb{R} ^{n-k|n-k} })\otimes\mathcal{O}_{ X \times \mathbb B^{n|n}})=(2n-2k)$$
		as a result,
		$$
		dim \,ker (\bar{D}(G_U \circ \Psi)\otimes 1_{\mathcal{O}_W})= (2m+2n)- (2n-2k)= (2m+2k)
		$$
		In addition, $i_W$ is immersion, thus we have
		$$dim( \bar{D}(i_W)\mathcal{T}_W(\Psi^{-1}(U)) =dim Der(\mathcal{O}_W)= (2m+2k).$$
		This completes the proof.	
	\end{proof}
	According to the diagram
	\ref{5.4.}
	$ ker (\bar{D}(G_U \circ \Psi)\otimes 1_{\mathcal{O}_W})$ is $ \Pi $-symmetric and its $\Pi$-symmetry is in the form
	$P_{X \times \mathbb B^{n|n}} \otimes 1_{\mathcal{O}_W}$,
	therefore
	$P_{X \times \mathbb B^{n|n}}|_{\mathcal T_W}$ is $ \Pi $-symmetry on 
	$\mathcal T_W(\Psi^{-1}(U))$
	for
	$U$'s
	are small enough. 
	\begin{proposition}
		$\mathcal T_{W}$
		is
		$\Pi$-symmetric.
	\end{proposition}
	\begin{proof}
		According to  \ref{5.12.}, for sufficiently small neighborhood $U_i$, there exists $\Pi$-symmetry as follows:
		$$P_{X \times \mathbb B^{n|n}}|_{\mathcal T_W}:\mathcal T_W(\Psi^{-1}(U_i)) \to \mathcal T_W(\Psi^{ -1}(U_i)).$$
		Now we construct a $\Pi$-symmetry on $W$ say $P_W$. \\To this end we need to construct$P_W(V):~\mathcal{T}(V)\to~\mathcal{T}(V)$ for any open subset $V$ of $X\times \mathbb B^{n|n}$. Set 
		$$ V= \bigcup U_i$$
		for $U_i$'s sufficiently small open subsets of $V$. Let $F\in \mathcal{T}_W(V)$. Set
		$$F_{U_i} =r_{U_i,V}(F)\,\,\,\,\, (U_i \subset V).$$
		Then one has
		$$r_{U_{ij}, U_i} (F_{U_i})=r_{U_{ij}, U_i}(r_{U_{i,V}}(F))=r_{U_{ji},V}(F)$$
		where
		$U_{ij}=U_i \cap U_j.$ On the other hand, we have
		\begin{align}
			& r_{U_{ij}, U_j} (F_{U_j} ) = r_{U_{ij}, U_j} ( r_{U_j, V} (F))=r_{U_{ji},V}(F)\\
			& \Rightarrow r_{U_{ij}, U_i} (F_{U_i}) = r_{U_{ji}, U_j} ( F_{U_j}).
		\end{align}
		Now we define $P_W(F)$ as follows:
		$$
		r_{U_j,V} ( P_W(F)):= P_{U_j} ( r_{U_j,V} (F))= P_{U_j} (F_{U_j})
		$$
		The vector field $P_W(F)$ is well defined if $P_{U_j}(F_{U_j})$ and $P_{U_i}(F_{U_i})$ have same restriction on $U_{ij}$. This is indeed the case. The following equalities complete the proof.
		\begin{align}
			& r_{U_{ij}, U_i} ( P_{U_i} (F_{U_i} ))= P_{U_{ij}} ( r_{U_{ij}, U_i} (F_{U_i} ))\\
			& r_{U_{ij}, U_j} ( P_{U_j} (F_{U_j}))= P_{U_{ij}} ( r_{U_{ij}, U_j} (F_{U_j}))
		\end{align}
	\end{proof}
	\section{Conclusion}
	Our goal in this article is to generalize Thom's transversality theorem in supergeometry which is one of the most important concepts in differential geometry. In this regard, we first generalized the concept of transversality in supergeometry and examined its stability property. In the remainder of the paper, we prove the genericity property through a generalization of Sard's theorem within the framework of supergeometry, genericity is established as the most important property for supertransversality in supergeometry 
	In the last section, we focus on a category of $\Pi$-symmetric supermanifolds and examine supertransversality 
	in this category.
	\appendix
	\section{Appendix}
	\textbf{Presheaf}.
	If $X$ is a topological space, then a presheaf is a correspondence, say $R$, that assigns to each open subset of $X$, say $U$, an algebraic structure $R(U)$. The correspondence $R$ must satisfy the following condition: if $V$ is an open subset of $X$ and $U$ is an open subset of $V$, then there exists a homomorphism $$r_{UV}: R(V) \to R(U)$$
	This homomorphism is called the restriction map.
	Another important condition is that if $U \subset V \subset W$ are open subsets, then we have
	$$r_{UW}=r_{UV} \circ r_{VW}$$
	where $\circ$ denotes the composition of homomorphisms.\\
	\textbf{Sheaf}.
	A sheaf is a presheaf say $R$ with the following condition: If $U$ is an open subset and $\{U_{\alpha}\}_\alpha$ is an open cover of $U$ and $\{f_{\alpha}\}$ is  a family of elements 	$f_\alpha \in R(U_\alpha)$ such that for any $ \alpha,\beta$ with nonempty intersection $U_\alpha \cap U_\beta$ one has
	$$f_\alpha|_{U_\alpha \cap U_\beta} = f_\beta|_{ U_\alpha \cap U_\beta}$$
	then there exists a unique element $f \in R(U)$ such that $f|_{U\alpha}= f_\alpha$.\\
	\textbf{Gluing}. Let $X_i$ be an open cover of $X$ and let $\mathcal{O}_i$ be a sheaf of rings on $\tilde X_i$. In addition assume $ f_{ij}: (X_{ji},\mathcal{O}_j|_{X{ji}})\to (X_{ij},\mathcal{O}_i|_{X{ij}})$ is an isomorphism of ringed spaces with $\tilde{ f_{ij}}=id|_{X_{ij}}$, where $X_{ij}=X_i \cap X_j$ for each $i,j.$ We call a ringed space $(X, \mathcal{O})$ is constructed by gluing ringed spaces $(X_i, \mathcal{O}_i)$ through $f_{ij}$ if there exists ringed space isomorphisms $$f_i: (X_i, \mathcal{O}|_{X_i})\to(X_i, \mathcal{O}_i|_{X_i})$$ with $\tilde{f}_i=id_{X_i}$ such that $f_{ij}=f_i\circ f_j^{-1}$. 
	\begin{theorem} \label{1.6.}
		Let $U^{m|n}$ be a superdomain. Take local coordinates $(t^i,e^j)$ , $ 1 \leq i \leq m $  and $ 1 \leq j \leq n$, for $U^{m|n}$. Suppose $X$ is a supermanifold. 
		$f_i$ and $g_j$ are even and odd elements of $\mathcal{O}_X(\tilde X)$, respectively. Then there exists a unique morphism $\psi:X\to U^{m|n}$ such that $\psi^*(t^i)= f_i$ and $\psi^*(e^j)=g_j$\cite{1}.
	\end{theorem}
	\textbf{Tangent space}.
	Let $X^{m|n}$ be a supermanifold and let $x \in \tilde X$. Then as in the classical case, a tangent vector to $X$ at $x$ is defined as a derivation of the stalk. If $(x_i,\xi_j)$ is a coordinate system for $X$ at the point $x$, then the tangent space to $X$ at $x$, denoted by $T_x(X)$ has $\{(\dfrac{\partial}{\partial x_i})_x,(\dfrac{\partial}{\partial \xi_j})_x \} $ as a basis and so is a super vector space of dimension $m|n$. If $\psi:X \to Y$ is a morphism of supermanifolds and $x \in \tilde{X}$ and $y=\tilde\psi(x) \in \tilde{Y}$, then tangent map of $\psi$ at $x$ is defined as follows:
	\begin{align*}
		d\psi_x: T_x(X) \to T_y(Y)\\
		\eta \mapsto \eta \circ \psi^*.
	\end{align*}
	In terms of coordinats $(x_i,\xi_j)$ for $X$ and $(y^a,\theta^b)$ for $Y$ with $\psi^*(y^a)=f_a$ and $\psi^*(\theta^b)=g_b,$ we define Jacobian of $\psi$ with a slight modification as  
	$$J\psi=
	\begin{pmatrix}
		\,\, \dfrac{\partial f_a}{\partial x_i}\,\,\,\ {-\dfrac{\partial f_a}{\partial \xi_j}}\\ \dfrac{\partial g_b}{\partial x_i}\,\,\,\,\,\,\,\ \dfrac{\partial g_b}{\partial \xi_j}	\\
	\end{pmatrix}
	.
	$$
	However the matrix of $d\psi_x$ with respect to coordinate basis of $T_xX$ and $T_yY$ is as follows:
	$$
	\begin{pmatrix}
		\,\, {\dfrac{\partial f_a}{\partial x_i}^\sim}(x)\,\,\,\ 0\,\,\\\ \,\,0\,\,\,\,\ {\dfrac{\partial g_b}{\partial \xi_j}}^\sim(x)	\\
	\end{pmatrix}.
	$$
	\begin{definition}(\textbf{Immersion}) \label{2.8.} 
		A morphism $\psi$ is called an immersion at $x$ if $d\psi_x$ is injective. A standard immersion is a morphism $\psi:U^{m|n} \to V^{m+r|n+s}$ with coordinates $(x^i,\xi^j)$ and $(x^i,y^a,\xi^j,\theta^b) $ for $U$ and $V$, respectively such that $\tilde{\psi}: x\mapsto (x,0)$ and
		$$\psi^*:x^i \mapsto x^i, \, \xi^j \mapsto  \xi^j, y^a ,\theta^b \mapsto 0$$
	\end{definition}
	\textbf{Embedding}. 
	Let $\iota:=(\tilde{\iota},\iota^*):X=(\tilde{X},\mathcal{O}_X) \to Y=(\tilde{Y},\mathcal{O}_Y) $ be a morphism of supermanifolds. We say that  $(X,\iota)$ is an embedded supermanifold if $\tilde{\iota}$ is an immersion and $\tilde{\iota}:\tilde{X}  \to \tilde{Y} $ is an homeomorphism onto its image.\\
	In particular, if $\iota  (X)\subset Y $ is a closed subset of $Y$ we will say that $(X,\iota)$ is a closed embedded supermanifold.\\
	In what follows, we will always deal with closed embedded supermanifolds. Remarkably, it is possible to show that a morphism  $\iota:X=(\tilde{X},\mathcal{O}_X) \to Y=(\tilde{Y},\mathcal{O}_Y)$ is an embedding if and only if the corresponding morphism $\iota^*:\mathcal{O}_Y \to \mathcal{O}_X $ is a surjective morphism of sheaves. Notice that, for example, given a supermanifold $Y=(\tilde{Y},\mathcal{O}_Y) $, one always has a natural closed embedding: the map $\iota:(\tilde{Y},\mathcal{O}_Y)_{red} \to (\tilde{Y},\mathcal{O}_Y) $, that embeds the reduced manifold underlying the supermanifold into the supermanifold itself \cite{14}.
	\begin{definition}(\textbf{Submersion}) \label{2.9.}
		A morphism $\psi$ is called a submersion at $x$ if $d\psi_x$ is surjective. A standard submersion is a morphism $\psi: V^{m+r|n+s} \to U^{m|n}$, with coordinates $(x^i,y^a,\xi^j,\theta^b) $ and  $(x^i,\xi^j)$ for $V$ and $U$, respectively such that $\tilde{\psi}: (x, v)\mapsto x $ and
		$$\psi^*:x^i \mapsto x^i, \, \xi^j \mapsto  \xi^j $$
	\end{definition}
	\begin{theorem} \label{2.10.}
		Standard immersions and submersions in definitions \ref{2.8.} and \ref{2.9.} are local models for immersion and submersion respectively.
	\end{theorem}
	\bibliographystyle{amsplain}
	
\end{document}